\documentclass[12pt,a4paper]{article}
\usepackage{amsmath}
\usepackage{graphicx}
\usepackage{amssymb}
\usepackage{hyperref}
\begin{document}
\title{Regular polygons}
\author{J. Mainik}
\date{November 15, 2026}
\maketitle
\begin{abstract}
The construction of regular polygons with a compass and straightedge 
is a problem that has interested mathematicians for a long time. 
The question of whether is it possible to construct the regular 17-gon,\linebreak
for example, remained unanswered for more than 2000 years. 
Gauss solved this problem in 1796. He proved later that it is possible to construct  
the regular polygons with $n=2^m n_1\cdots n_{\ell}$ sides, where $n_1,\cdots, n_{\ell}$ 
are different Fermat primes. 

This paper presents a new approach for solving this task.
Among other things we analyze in detail the case of $n=65537$. 
J. G. Hermes announced in 1894 that he had a full description of the construction 
of the 65537-gon. But his text was too extensive and was never published. 
We show exactly and without gaps how the regular 65537-gon can be constructed.
\end{abstract}

\section{Introduction}

The construction of regular polygons with a compass and straightedge 
is a  well-known task which has engaged mathematicians for a long time.
Particularly interesting was the question of whether is it possible to construct
the regular 17-gon.
This task, which is immediately understandable and 
can be formulated so simply, 
could not be solved for more than 2000 years. C. F. Gauss solved this problem in 1796,
[1,2].
He proved later that it is possible to construct with a compass and straightedge 
the regular polygons with $n=2^m n_1\cdots n_{\ell}$ sides, where $n_1,\cdots, n_{\ell}$ 
are different Fermat primes $n_k=2^{2^{\nu_k}}+1$.
The numbers of the form $n=2^{2^{\nu}}+1$ are known as Fermat numbers.
Fermat thought that these numbers are all prime numbers.
However, that is not true: the only currently known Fermat primes are
$3$, $5$, $17$, $257$ and $65537$.\newpage

Gauss came to the solution by his research in the area of number theory. 
He made use of the fact that for a prime number $n$ it is possible to find a primitive root.
The number $g$ is a primitive root for $n$, if 
the remainders $r_k=\mathrm{rest}(g^k,n)$ of $g^k$ divided by $n$ for $1\le k\le n-1$
give all numbers $1,2,3,\cdots,n-1$. 

For the prime number $n=17$ Gauss considered the ordered set
\[
e^{\frac{2\pi i}{n}r_k}, 1\le k\le 16,
\]
and divided it into two smaller parts so that the elements were added 
to the parts alternately. The smaller parts he split again in the same way.
Gauss noticed that this gives rise to sets for which the values can be calculated
with the help of quadratic equations and 
therefore constructed with a compass and straightedge. 
This is how he constructed $\cos(2\pi/17)$ and the regular $17$-gon.

Gauss also knew that only these regular polygons can be constructed but did not prove it.
P. Wantzel completed the result of Gauss and proved it in 1837,
[6].

It is clear that a regular polygon with $n=2^m n_1\cdots n_{\ell}$ sides,
 where $n_1,\cdots, n_{\ell}$ are different Fermat primes,
can be constructed if and only if every $n_k$-gon, $1\le k\le \ell$, can be constructed.
Therefore, we consider in the following only regular $n$-gons for Fermat primes.
We present an approach that differs from the method of Gauss. 
In particular, we show  how the regular 65537-gon can be constructed.

There is a nearly infinite list of publications on this subject.
We show only a few of them,
[1-16], 
but this list confirms very well that work is being carried out on this task continuously.
This is so because the search is still going on for shorter and prettier 
constructions of the regular polygons.

One of the first geometric constructions of the 17-gon was presented
by M. G. Paucker in 1819 and published in 1822,
[3].
A simpler construction of the 17-gon was found by H. W. Richmond in 1893,
[7].

Paucker proposed the first description of the construction for the 257-gon\linebreak
 in 1822,
[3].
In 1825 J. Erchinger proposed an another construction,
and this construction was discussed by Gauss himself,
[4].
In 1832 the same construction was once again presented by F. J. Richelot, 
[5]. 
D. W. DeTemple in the year 1991, M. Trott in 1995 and C. Gottlieb in 1999 
published other constructions of the 257-gon,
[11-13].

It took a long time until a construction of the regular 65537-gon was presented.
In the year 
1894, nearly $100$ years after the publication of Gauss, \linebreak
J. Hermes announced that he had finished his work and had a full 
construction of the regular 65537-gon. This was the result of 10 years of work,
but the paper of Hermes was too extensive and was never published. This paper 
is kept in the library of the University of Goettingen.
Hermes was able to publish only a 17 page summary in 1895,
[6].

It is clear that the construction of regular polygons with a compass and straightedge 
has not the highest practical relevance. 
But the interest of mathematicians and non-mathematicians for this task is still present.
That has to do with the elegance of the regular $n$-gons and the simplicity
of the question, but even more with the need for solutions of such understandable tasks.

The author hopes that the new approach to constructing the regular polygons
and the presented construction for the 65537-gon will find interest and understanding. 
\section{Denominations and remarks}

In the following, $n=2^{2^ {\nu}}+1$ is a Fermat prime. For $z=e^{2\pi i/n}$ we have
\begin{eqnarray*}\label{Gl1}
z^n-1=0
\end{eqnarray*}
and the the values $z^k$ for $k=0,1,\cdots ,n-1$ are the solutions of this equation.
The points $z^k=e^{k\cdot 2\pi i/n}$, $0\le k\le n-1$, are lying on the circle with 
the radius \nolinebreak $1$ 
and represent the vertices of the regular $n$-gon. The radius 
of this circle equals $1$, and therefore we have $z^0=1$.
It applies
\begin{eqnarray*}\label{SummexHochk}
z^n-1=(z-1)(z^{n-1}+z^{n-2}+\cdots +z+1)\nonumber
\end{eqnarray*}
and the values $z^k=e^{k\cdot 2\pi i/n}$, $1\le k\le n-1$, which we need also, are the
solutions of the equation
\begin{eqnarray}\label{SummexHochk}
z^{n-1}+z^{n-2}+\cdots +z+1=0\;.
\end{eqnarray}
In order to construct the regular $n$-gon we have thus to determine the solutions of this 
equation, and these solutions should be constructed with a compass and straightedge.

In the following we always designate the terms $z^k=e^{k\cdot 2\pi i/n}$, $1\le k\le n-1$, 
as elements.

In this context we summarize a few simple facts. The solutions\linebreak
 $x_{1,2}=-p/2 \pm \sqrt{p^2/4-q}$ of the quadratic equation $x^2+px+q=0$ 
can be constructed with a compass and straightedge, if they are real values and we have 
the segments with the length $p$ and $q$. For this construction
we need in addition to these segments the segment of length $1$, but this 
is the radius of the circle on which the vertices of the 
$n$-gon are positioned.
Here a remark about the square and the square root may be worthwhile.
The value $p^2$ can be constructed using the equality $p^2/p=p/1$
and the value $\sqrt D$ can be constructed as follows: Construct the semicircle with diameter
$D+1$ and the perpendicular to this diameter at the point that splits the
diameter into segments of lengths $D$ and 1. The perpendicular line
from this point to the semicircle has length $\sqrt D$.

By constructing with a compass and straightedge we are, in a certain sense,
limited to solutions of quadratic equations. It is known, for example,
that the solution of the simple cubic equation
\[
x^3-3x-1=0
\]
that appears by trisection of the 60-degree angle and even the root $\sqrt[3]{2}$
of the equation $x^3=2$ that we need for doubling the cube cannot 
be constructed with a compass and straightedge.
In order to determine the solution of a more complicated equation we have
to come repeatedly to quadratic equations and construct the solutions of them.

For the quadratic equation with solutions $x_1$ and $x_2$ we have
\[
(x-x_1)(x-x_2)=x^2-(x_1+x_2)\cdot x+(x_1\cdot x_2)=0
\]
and to determine this quadratic equation, we must be able to determine the sum
$x_1+x_2$ and the product $x_1\cdot x_2$.

In the following we will always split already known, complicated sums so that all terms 
of the complicated sum will be disjunct distributed between two parts.
The sum of the parts is then automatically equal to the known value
and we always just have to determine the product of the parts.
For the fairly complicated equation
(\ref{SummexHochk}),
for example, we must therefore several times split the sum $S$ 
\begin{equation}\label{S=}
S=z^{n-1}+z^{n-2}+\cdots +z
\end{equation}
into smaller parts in such a way that we can determine the products of these parts.
Then we can construct the parts with a compass and straightedge.

In the following we will always have parts with real values, and, 
to assign the real roots of the quadratic equation to the geometric objects correctly,
we need to know which of the values is to calculate with the sign "+" 
and which is to calculate with the sign "$-$" before the discriminant.
That is, we need to know which of them is greater.

\section{Parts of the sum $S$ and invariant sets}

If $F_1$ and $F_2$ are disjunct parts of $S$, it holds $F_1+F_2=S$,
and we have to choose  $F_1$ and $F_2$ in a way that allows us 
to determine the product $F_1\cdot F_2$.
We know the value $S=-1$, but we don't know any part of $S$.
It is therefore natural to choose the parts $F_1$ and $F_2$ so 
that in the product $F_1\cdot F_2$
each term $z^k$, $1\le k\le n-1$, of $S$ appears an equal amount of times.
If each $z^k$ appears $\mu$ times we say that we have an even coverage of height $\mu$.
In this case it holds $F_1\cdot F_2=\mu\cdot S+\delta$. 
The constant summand $\delta$ here will be the number
of elements $z^k$ in $F_1$ for which the inverse elements $z^{n-k}$ belong to $F_2$.\linebreak
In this way we can obtain the necessary information to determine the relevant
quadratic equation for the parts $F_1$ and $F_2\,$.

From  
(\ref{Gl1})
we see that $z^{n+m}=z^{m}$, and therefore we get 
\begin{eqnarray*}
z^k\cdot S &=&z+z^2+\cdots z^{k-1}+z^{k+1}+\cdots+z^{n-1}+1
\end{eqnarray*}
so that
\begin{eqnarray}\label{z^kS}
z^k+z^k\cdot S=S+1\;.
\end{eqnarray}
If we add up in
(\ref{z^kS})
over all values $z^k$ that belong to $F_1$, we obtain
\begin{eqnarray*}
F_1+F_1\cdot S=\mu(F_1)\cdot S+\mu(F_1)\;,
\end{eqnarray*}
where $\mu(F_1)$ denotes the number of elements of $F_1$.
Because of $S=F_1+F_2$, we have then
\begin{eqnarray}\label{F+FF}
(F_1+F_1^2)+F_1\cdot F_2=\mu(F_1)\cdot S +\mu(F_1)\;.
\end{eqnarray}
From 
(\ref{F+FF})
we see that $F_1\cdot F_2$ provides an even coverage of $S$ 
if and only if $F_1+F_1^2$ provides an even coverage of $S$. 
It is natural, in this context, 
to take a closer look at the squares of the elements of the part $F_1\,$.
\newline

{\bf Proposition 1.} Let $n>5$ and let $F_1$ be so that $F_1+F_1^2$ provides
an even coverage of $S$. Then with $z^k\in F_1$ it also holds $z^{2k}\in F_1\,$.

{\bf Proof.} At first we notice that $z_1^2\ne z_2^2$ if
$z_1$ and $z_2$ are different elements, $z_1\ne z_2$. This property is obvious.

Let $z_1\in F_1$ and $z_1^2\not\in F_1$. The element $z_1^2$ will then
obviously have an odd coverage with $F_1+F_1^2$ : once as the element
$z_1^2$ in $F_1^2$ and in addition possibly an even number of times 
by products $2z_{\ell}\cdot z_m$ of different elements $z_{\ell}$ and $z_m$ of $F_1$
by the calculation of $F_1^2$. It follows immediately that for
each element $z_k$ of $F_1$ the square $z_k^2$ cannot belong to $F_1$.
It has to be so because otherwise the element $z_k^2$, $z_k^2\in F_1$, will be
covered with $F_1+F_1^2$ an even number of times.
As $F_1+F_1^2$ provides an even coverage, all elements of $S$ 
must have this odd coverage. This means that in $F_1+F_1^2$ each element 
of $F_2$ appears exactly once as $z_k^2$ for $z_k\in F_1$.
Thus we see that $F_1$ and $F_2$ must have the same number $\frac{n-1}{2}$ 
of elements and the function $f(z)=z^2$ provides a one-to-one mapping $F_1$ to $F_2$. 

We see in
(\ref{F+FF})
at once that together with $F_1+F_1^2$ the product $F_1\cdot F_2$
provides an odd coverage of $S$.
The product $F_1\cdot F_2$ has $(\frac{n-1}{2})^2-\delta$
elements, where $\delta$ is the number of elements $z^k$ in $F_1$ for which
the inverse elements $z^{n-k}$ belong to $F_2$.
It holds here obviously $0\le \delta\le \frac{n-1}{2}$.
As $F_1\cdot F_2$ provides an even coverage of $S$, $\:n-1$ is a divisor 
of $(\frac{n-1}{2})^2-\delta$ and this is possible only if $\delta=0$.
Thus we see that the coverage $\mu$ for $F_1\cdot F_2$ is $\mu=\frac{n-1}{4}$.
We come to a contradiction as $\frac{n-1}{4}=2^{2^{\nu}-2}$ is an even number.
$\blacksquare$
 
{\bf Corollary.}  As a corollary of the proved proposition one sees that with each element
$z^k$ of $F_1$ all elements $z^{2k}, z^{4k}, z^{8k},\cdots$ also belong to $F_1$.
\newline

{\bf Definition. } An invariant set is a set of elements, formed in such a way that, 
in addition to  a starting element $z^s$, all the distinct elements, 
obtained through the calculation of squares, also belong to the set.
An invariant set therefore consists
of the elements $z^s, z^{2s}, z^{4s}, \dots, z^{2^ks}$, 
where $k\ge 1$ is the smallest integer such that $z^{2^{k+1}s}$ 
returns to an already present value: $z^{2^{k+1}s} = z^{2^js}$, $0 \le j \le k$.
The order of elements of an invariant set, where the next element
is the square of the previous one, we denote a natural order.\newline

{\bf Proposition 2.} If $n\ge 5$, then every invariant set consists of $n_1=2^{\nu+1}$ 
different elements. By natural order the $m$-th element, $1\le m\le 2^{\nu}$, 
of the invariant set is inverse to the  $(2^{\nu}+m)$-th element of the set.

{\bf Proof. }The invariant set, starting with the element $z^s$, $1\le s\le n$, 
is formed by calculating the corresponding squares. In this way, we obtain the values
$z^{2^k\cdot s}$, $k = 0, 1, 2, \dots$, and must take into account that $z^n = 1$.
For the degrees of these elements we then have the values
$2^k\cdot s$, $k=0, 1, 2, \dots$,  and must calculate {\it modulo }$n$.

For $0\le k\le 2^{\nu}$ the values $2^k$ increase from 1 to $2^{2^{\nu}}=n-1$.
For $k=2^{\nu}+k_1$, $k_1\ge 1$, it holds
\[
2^{2^{\nu}+k_1}=2^{2^{\nu}}\cdot 2^{k_1}=(n-1)\cdot 2^{k_1}\equiv n-2^{k_1}\nolinebreak \pmod{n}\;,
\]
and we see that for $2^{\nu}+1\le k\le 2^{\nu+1}$ 
the values $\mathrm{rest}(2^k,n)$ decrease 
from $n-2$ to $1$. The value $\mathrm{rest}(2^k,n)$, $0\le k\le 2^{\nu+1}$, 
is equal to 1 only for $k=0$ and $k=2^{\nu+1}$, and for $1\le k< 2^{\nu+1}$
it then holds $2^k\not\equiv 1\pmod{n}$.

For the degrees $2^{k_1}\cdot s$ and  $2^{k_2}\cdot s$, where $0\le k_1<k_2< 2^{\nu+1}$,
the following holds:
\begin{eqnarray*}
&&2^{k_2}\cdot s -2^{k_1}\cdot s=2^{k_1}\cdot(2^{k_2-k_1}-1)\cdot s\not\equiv 0\pmod{n}\;,
\end{eqnarray*}
and the corresponding elements of the invariant set are different.\linebreak
For $k=2^{\nu+1}$ we get $ 2^{2^{\nu+1}}\cdot s=(n-1)^2\cdot s\equiv s\pmod{n}$, 
and we see that we come back to the initial element $z^s$.
The invariant set thus consists of $n_1=2^{\nu+1}$ different elements.

For the sum of the degrees of the $k$-th and the $(2^{\nu}+k)$-th element in the
invariant set we have
\[
2^{k-1}\cdot s+ 2^{2^{\nu}+k-1}\cdot s= 2^{k-1}\cdot (2^{2^{\nu}}+1)\cdot s=
2^{k-1}\cdot n\cdot s\equiv 0\pmod{n}\;,
\]
and that obviously means that these elements are inverse.
$\blacksquare$
\newline

If we were to calculate the squares again and again, they would repeatedly 
run through the values of the invariant set in the same order after returning 
to the starting element.
For the author, this property was the reason for the name "invariant set."

It is clear that the invariant sets cannot overlap. From the overlapping
point they will coincide, and, due to the coming back to the starting element,
they should be the same.
This means that the $2^{2^\nu}$ elements of $S$ are distributed between 
$2^{2^{\nu}}/2^{\nu+1}=2^{2^{\nu}-(\nu+1)}$ disjunct invariant sets.
In the following, we denote the number of the invariant sets by $n_2$,
$n_2=2^{2^{\nu}-(\nu+1)}$.

Each invariant set consists of $n_1/2$ pairs of inverse elements,
and it is natural to unite the inverse elements in pairs and not separate them.
There are two reasons for this. At first, when working with pairs,
we are dealing with real values 
\[
p_k=z^k+z^{n-k}=2 \cos\frac{2k\pi}{n}\; ,
\]
and we are particularly interested in the value $p_1=2 \cos\frac{2\pi}{n}$.
If we determine this value, we can simply construct the polygon completely.
The second reason for the work with pairs is that we will have half as many values.

We have assumed above that for the pair $p_k=z^k+z^{n-k}$ 
the number $k$ is chosen as $k=\mathrm{min}(k,n-k)$ and therefore 
the numbers $k$ of the pairs are in the range  $1\le k\le\frac{n-1}{2}$.
So we get a unique numbering of the pairs. 
One can also represent the pair $p_k$ as $p_k=z^k+z^{-k}$,
and this form can be more appropriate for some calculations.
The natural order for the pairs is such that the pair $p_k$ is followed 
by the pair $p_{\mathrm{min}(2k,n-2k)}$. 
We have to convert the element number $2k$ into a pair number.

To work with pairs, we need to know how we can multiply them.
Let $p_k=z^k+z^{-k}$, $1\le k \le\frac{n-1}{2}$,
and $p_m=z^m+z^{-m}$, $1\le m \le\frac{n-1}{2}$, be two pairs.
Then we have
\[
p_k\cdot p_m=(z^k+z^{-k})(z^m+z^{-m})=z^{k+m}+z^{-(k+m)}+ z^{k-m}+z^{m-k}\; ,
\]
and we see that the product of two different pairs is the sum of two pairs
and the square of a pair is the sum of a pair and the constant $2\,$:
\begin{equation}\label{pkpm}
p_k\cdot p_m=\begin{cases}
p_{\mathrm{min}(k+m,n-(k+m))}+p_{|k-m|}, & \texttt{if }p_k\neq p_m\\
p_{\mathrm{min}(2k, n-2k)}+2, & \texttt{if }p_k=p_m
\end{cases}
\end{equation}

We can consider the elements $z^k$ and the pairs $p_k=z^k+z^{-k}$
for any integer $k$ with $k\not\equiv 0\pmod{n}$. In this case  
\[
z^{k_1}=z^{k_2}\Leftrightarrow k_2\equiv k_1 \pmod{n}
\]
and
\[
p_{k_1}=p_{k_2} \Leftrightarrow k_2\equiv k_1 \textrm{ or }\;\; k_2\equiv -k_1\pmod{n}
\]
obviously applies. This is appropriate when we are interested 
in the value, and not in the exact number of the element or the pair.
This allows us to omit the calculation {\it modulo }$n$ and, when working with pairs,
also the conversion of the element number into a pair number.
For the product of pairs $p_k$ and $p_m$ with $k\not\equiv 0\pmod{n}$ 
and $m\not\equiv 0\pmod{n}$ we then have
\begin{equation*}
p_k\cdot p_m=\begin{cases}
p_{k+m}+p_{|k-m|}, & \texttt{if }p_k\neq p_m\\
p_{2k}+2, & \texttt{if }p_k=p_m
\end{cases}
\end{equation*}
In the following we will indicate which numbers are being used, if this is not obvious.

We further denote the different invariant sets by $E_k$, $1\le k\le n_2$.
The notation $E_k$ is intended to indicate that these are sets of elements. \linebreak
We denote the corresponding sets of the $n_1$ degrees of the elements from $E_k$
by $\widehat{E}_k$ and the corresponding sets of the $n_1/2$ numbers 
of the pairs from $E_k$ by $\Bar{E}_k$. 
So we have $E_k=\{z^k|k \in \widehat{E}_k\}=\{p_k|k \in \Bar{E}_k\}$.
The invariant sets $E_k$ form a disjunct partition 
of the elements $z^1,z^2,\cdots,z^{n-1}$,
the sets $\widehat{E}_k$ form a disjunct partition of the numbers $1,2,\cdots,n-1$,
and the sets $\Bar{E}_k$ form a disjunct partition of the numbers $1,2,\cdots,(n-1)/2$.
To emphasize the connections, we also call the sets $\widehat{E}_k$ and $\Bar{E}_k$ 
invariant sets.

For the invariant set $E_k$ we denote the sum of the corresponding elements by $S_k$ 
and call it a block or an invariant block.
The name "block" is intended to indicate that the elements are grouped together 
in a particular way, and the notation $S_k$ is meant to show that this is a sum.
So it applies \linebreak
 $S_k=\sum_{k\in \widehat{E}_k}z^k=\sum_{k\in \Bar{E}_k}p_k$.
Since $E_k$ consists of inverse pairs, $S_k$ is a real number
$S_k=\sum_{k\in \Bar{E}_k}\cos(2\pi k/n)$.

Construction of the regular $3$-gon is trivial, and construction of the regular 
$5$-gon is also a simple planimetric task. However, for the sake of completeness, 
we want to consider these cases and show how the approach works.

For $n=3$, there is only one invariant set $E_1$  with one pair $\bar{E}_1=\{1\}$,
 and we immediately have $p_1=2\cos(2\pi/3)=-1$.
Using $\cos(2\pi/3)=-1/2$, we can construct the angle $2\pi/3$ and, consequently, 
the regular triangle.

For $n=5$, there is only one invariant set $E_1$  with two pairs: $\bar{E}_1=\{1,2\}$.
Here we have $p_1+p_2=-1$ and we can split $p_1+p_2$ only into $p_1$ and $p_2$.
Using 
(\ref{pkpm})
we obtain $p_1\cdot p_2=p_1+p_2=-1$, and the quadratic equation for $p_1$ and $p_2$ is
\[
x^2+x-1=0\;.
\] 
Since $p_1>p_2$, we have $p_1=(-1+\sqrt{5})/2$.
We can construct  $p_1$ and, using $\cos(2\pi/5)=p_1/2$,  construct the angle $2\pi/5$ 
and the regular pentagon.

\section{Building invariant sets}

Let $n\ge 17$ be a Fermat prime, and let $g$ be a primitive root for $n$.
For any number $s$, $1\le s\le n-1$, the following holds for $1\le j\le n-1$ :
\[
s\cdot g^j\equiv s\Leftrightarrow s\cdot(g^j -1)\equiv 0\Leftrightarrow g^j\equiv 1
\Leftrightarrow j=n-1 \pmod{n}
\]
and $\mathrm{rest}(s\cdot g^j,n)$ returns to the starting value $s$
for the first time at $j=n-1$.
For $0\le j_1<j_2< n-1$ we then have
\[
s\cdot g^{j_2}-s\cdot g^{j_1}= s\cdot g^{j_1}\cdot (g^{j_2-j_1}-1)\not \equiv 0\pmod{n}\;.
\]
Thus we see that the values $\mathrm{rest}(s\cdot g^j,n)$, $0\le j<n-1$, 
are all different, cover 
all values $1,\cdots,n-1$, and therefore reach into all invariant sets $\widehat{E}_k$.
We want to analyze how the invariant sets $\widehat{E}_k$ themselves 
are traversed in the process.

If $s$ belongs to $\widehat{E}_k$, then by using $\mathrm{rest}(s\cdot 2^m,n)$, $0\le m<n_1$,
we can determine all values from $\widehat{E}_k$. 
Therefore, $\mathrm{rest}(s\cdot g^j,n)$ returns to $\widehat{E}_k$ from $s$
if there exists a number  $m$, $0\le m<n_1$, such that it holds
\[
s\cdot g^j\equiv s\cdot 2^m \Leftrightarrow s\cdot (g^j-2^m)\equiv 0 \Leftrightarrow g^j\equiv 2^m \pmod{n}.
\]
In the condition $g^j\equiv 2^m \pmod{n}$ the initial value $s$ 
and the number $k$ of the invariant set $\widehat{E}_k$ play no role, 
and this means that, starting from any element in any invariant set,
we come back to the same invariant set in $j$ steps.

Let $j_0$ be the smallest positive integer such that for the corresponding
$m_0$ the congruence $g^{j_0}\equiv 2^{m_0}\pmod{n}$ holds, 
and let $\widehat{E}_{k},\cdots,\widehat{E}_{k+j_0-1}$  be the sets 
that were passed through during these $j_0$ steps.
Since $j_0$ is the smallest value, the sets $\widehat{E}_{k},\cdots,\widehat{E}_{k+j_0-1}$
are obviously different, and therefore $j_0\le n_2$ must hold. But we have
\[
s\cdot g^{j_0+j}= s\cdot g^{j_0}\cdot g^j\equiv s\cdot 2^{m_0} \cdot g^j\equiv 
s\cdot g^j\cdot 2^{m_0} \pmod{n}
\] 
and we see that $\mathrm{rest}(s\cdot g^{j_0+j},n)$  and 
$\mathrm{rest}(s\cdot g^{j},n)$  belong to the same invariant 
set $\widehat{E}_k$.
This means that with $\mathrm{rest}(s\cdot g^{j},n)$ for $j=1,2,\cdots$ 
the $j_0$ different sets $\widehat{E}_{k},\cdots,\widehat{E}_{k+j_0-1}$ 
are traversed repeatedly in the same order.
But as $\mathrm{rest}(s\cdot g^{j},n)$  appears in all invariant sets,
it must be $j_0\ge n_2$. Thus, $j_0=n_2$.

The sets $\widehat{E}_k=\{r_{k,1},r_{k,2},\cdots,r_{k,n_1}\}$, $1\le k\le n_2$,
with natural ordering can thus be calculated as follows.
With $r_{k,1}=\mathrm{rest}(g^{k-1},n)$, $1\le k\le n_2$, we determine the starting value in
$\widehat{E}_k$, and $r_{k,j}$, $2\le j\le n_1$, can be calculated using doublings:
$r_{k,j}=\mathrm{rest}(2r_{k,j-1},n)$, $2\le j\le n_1$.

The sets $\Bar{E}_k=\{s_{k,1},s_{k,2}\cdots,s_{k,n_1/2}\}$, $1\le k\le n_2$,
with natural ordering can be calculated analogously:
\begin{eqnarray*}\label{p_Number}
s_{k,1}&=&\mathrm{min}(\mathrm{rest}(g^{k-1},n),n-\mathrm{rest}(g^{k-1},n)),\; 1\le k\le n_2\; ,\nonumber\\
s_{k,j}&=&\mathrm{min}(\mathrm{rest}(2s_{k,j-1},n),n-\mathrm{rest}(2s_{k,j-1},n)),\; 
2\le j\le n_1/2\; .
\end{eqnarray*}

In this way we can also build the invariant sets $E_k$ for $k>n_2$. 
Doing so, we obviously obtain a simple cyclic property: $E_{n_2+1}=E_1,E_{n_2+2}=E_2,\cdots$.

For 17, 257 and 65537 we can use 3 as a primitive root.
But to avoid excluding any possibly existing larger Fermat prime, 
we built  the invariant sets using the not precisely known primitive root $g$.

In the case $n=17$, there are only 2 invariant sets, each with $4$ pairs:
$\Bar{E}_1=\{1,2,4,8\}$, $\Bar{E}_2=\{3,6,5,7\}$.

In the case $n=257$, we have $16$ sets  $\Bar{E}_k$, each with $8$ numbers, which we can 
easily calculate and use. We are showing only $4$ of them in detail:
\begin{eqnarray*}
\Bar{E}_1 &=& \{1,\;2,\;4,\;8,\;16,\;32\;,64,\;128\}\; , \\ 
\Bar{E}_2 &=& \{3,\;6,\;12,\;24,\;48,\;96,\;65,\;127\}\; ,\\ 
\Bar{E}_3 &=& \{9,\;18\,;36,\;72,\;113,\;31,\;62,\;124\}\; ,\\
\cdots\\
\Bar{E}_{16} &=& \{83,\;91\,;75,\;107,\;43,\;86,\;85,\;87\}\; .
\end{eqnarray*}
In the case $n=65537$, we have $2048$ sets $\Bar{E}_k$, each with $16$ numbers, 
and they cannot be represented so simply.

\section{Regular 17-gon}

The results of Propositions 1 and 2 are sufficient to construct the regular 17-gon.
Here we have only two invariant sets  $E_1$ and $E_2$, and it applies:
$S_1=p_1+p_2+p_4+p_8$ and $S_2=p_3+p_6+p_5+p_7$.
The natural order of pairs yields: $p_1^2=p_2+2, p_2^2=p_4+2,p_4^2=p_8+2 $.
It holds also $p_8^2=p_1+2$ so that we have a certain circular property 
for the pairs in $S_1$. The same is also true for the pairs in $S_2$.\pagebreak

According to Proposition 1, $S$ must be divided into $S_1$ and $S_2$.
We have the sum $S_1+S_2=S=-1$, and only the product $S_1\cdot S_2$ must be calculated.
We obtain via direct calculation
\begin{eqnarray*}
S_1\cdot S_2 =  4(p_1+p_2+p_3+p_4+p_5+p_6+p_7+p_8)\;=\;4S\;=\;-4
\end{eqnarray*}
and the corresponding quadratic equation is thus
\begin{eqnarray}\label{QuadratGl1}
x^2+x-4=0.
\end{eqnarray}
It holds here $S_1>S_2$, and we obtain therefore
\[
S_1=\frac{-1+\sqrt{17}}{2}, \;\;\;
S_2=\frac{-1-\sqrt{17}}{2}\;.
\]
These values can be constructed with a compass and straightedge.

We want to get the pair $p_1$, and hence it is natural to split the 
sum $S_1=p_1+p_2+p_4+p_8$ into disjunct parts $P_1$ and $P_2$,
so that $p_1$ belongs\linebreak
to $P_1$. The sum of the corresponding parts $P_1$ and $P_2$ is known: $P_1+P_2=S_1$.
As we know only $S_1$ and $S_2$ ($S$ is also the sum of $S_1$ and $S_2$),
we have to choose the parts $P_1$ and $P_2$ in a way that allows us to calculate 
the product $P_1\cdot P_2$ on the basis of $S_1$ and $S_2$. 
For this purpose, even coverings of $S_1$ and $S_2$ can be used.
There are not too many possibilities
to split $S_1$, and we can check them individually.
The splitting of $S_1$ where one part has 1 pair and the other part has 3 pairs
does not fit: We then get 6 pairs in the products, which cannot be
distributed evenly between $S_1$ and $S_2$ as they have 4 pairs each.
If we consider the parts with 2 pairs in each of them, we see simply that
the splittings into $P_1=p_1+p_2$ and $P_2=p_4+p_8$, or into $P_1=p_1+p_8$ 
and $P_2=p_2+p_4$ do not fit. Thus, it can be only the splitting of 
$S_1$ into $P_1=p_1+p_4$ and $P_2= p_2+p_8$.

If we calculate $P_1\cdot P_2$ for this split directly, we gain
\begin{eqnarray*}
P_1\cdot P_2 = (p_1+p_4)(p_2+p_8)=p_1+p_3+p_7+p_8+p_2+p_6+p_4+p_5=S=-1.
\end{eqnarray*}
Thus, the appropriate quadratic equation is here
\begin{eqnarray*}
x^2-S_1\cdot x-1=0
\end{eqnarray*}
and, as $P_1>P_2$, we obtain
\begin{eqnarray*}
P_1=\frac{S_1+\sqrt{S_1^2+4}}{2},\;\;\;
P_2=\frac{S_1-\sqrt{S_1^2+4}}{2}\;.
\end{eqnarray*}
These values can be constructed.
It can also be noted here that $S_1$ is a solution of the equation 
(\ref{QuadratGl1}), 
so that
\begin{equation}\label{S^2}
S_1^2=4-S_1
\end{equation}\nopagebreak
holds, and the values $P_1$ and $P_2$ can be expressed without $S_1^2$. \pagebreak

Here it is possible to calculate the product $P_1\cdot P_2$ 
for the split of $S_1$ with the help of squares. We obtain
\begin{eqnarray*}
&&\!\!\!\!\!\!\!\! 2P_1\cdot P_2 \; = \; S_1^2-(P_1^2+P_2^2)= S_1^2-((p_1+p_4)^2+(p_2+p_8)^2)\\
&&\!\!\!\;\!\!\!\! =\, S_1^2-(p_1^2+p_2^2+p_4^2+p_8^2)-2(p_3+p_5+p_6+p_7)
= S_1^2-S_1-8-2S_2\; .
\end{eqnarray*}
The value $(S_1^2-S_1-8-2S_2)/2$ can be constructed, and we certainly 
could determine the corresponding quadratic equation using this value.
But if we take 
(\ref{S^2})
into account, we obtain the well-known value
\[
P_1\cdot P_2=\frac{S_1^2-S_1-8-2S_2}{2}=\frac{4-6}{2}=-1\; .
\]
In any case we can notice here that it is appropriate
to split the invariant set by picking the pairs in alternating order.
Then we can easier calculate the squares of the relevant parts.

Now we have $p_1+p_4=P_1$ and want to split it into $p_1$ and $p_4$. 
If we calculate the product $p_1\cdot p_4$ directly, we obtain
\[
p_1\cdot p_4=p_3+p_5\; .
\]
The value $p_3+p_5$ presents a part of $S_2$, and
to determine this value we have to split $S_2$.
It can also be seen that the parts for $S_2=p_3+p_6+p_5+p_7$
should be $p_3+p_5$ and $p_6+p_7$, and for these parts the pairs
should be again chosen in alternating order.
This is certainly not difficult, but we want to get by even without that.
We want to use the squares. It applies obviously
\[
2p_1\cdot p_4=(p_1+p_4)^2-(p_1^2+p_4^2)=P_1^2-(p_2+p_8+4)=P_1^2-P_2-4
\] 
and therefore
\[
p_1\cdot p_4=\frac{P_1^2-P_2-4}{2}\; .
\]
Hence the quadratic equation for $p_1$ and $p_4$ is
\[
x^2-P_1\cdot x+\frac{P_1^2-P_2-4}{2}=0\;,
\]
and because of $p_1>p_4$ we obtain
\[
p_1=\frac{P_1}{2}+\sqrt{\frac{P_1^2}{4}-\frac{P_1^2-P_2-4}{2}}=
\frac{P_1+\sqrt{2P_2-P_1^2+8}}{2}\; .
\]
We can construct with a compass and straightedge the value 
$p_1=2\cos\frac{2\pi}{17}$ and therefore the complete regular 17-gon.
 
\section{Multiplication of invariant blocks}

Next we analyze what we will get as a product of invariant blocks. 
This is a natural question. Splitting $S$ into invariant blocks and
calculating the product of the parts,  we come to the  product of invariant blocks.\newline

{\bf Proposition 3. } The product of two different invariant blocks is the sum of
 $2^{\nu+1}$ invariant blocks.

{\bf Proof. }Let $S_k$ and $S_{\ell}$ be two different invariant blocks.
When calculating the product $S_k\cdot S_{\ell}$, we must multiply the $n_1$ 
elements $z^{r_{k,j_1}}$, $1\le j_1\le n_1$, from $E_k$ by the $n_1$ elements 
$z^{r_{{\ell},j_2}}$, $1\le j_2\le n_1$ from $E_{\ell}$ and obtain $n_1^2$ terms 
$z^{r_{k,j_1}+r_{{\ell},j_2}}$, $1\le j_1\le n_1$, $1\le j_2\le n_1$.
We omit the calculation {\it modulo }$n$ here because we are interested in the values
and not in the exact indices of the products. Instead of $r_{k,j_1}$, $1\le j_1\le n_1$,
we can obviously use the degrees $r_{k,1},2r_{k,1},\cdots,2^{n_1-1}r_{k,1}$,
and instead of $r_{{\ell},j_2}$, $1\le j_2\le n_1$, we can use the degrees
$r_{{\ell},1},2r_{{\ell},1},\cdots,2^{n_1-1}r_{{\ell},1}$.
These are suitable degrees because they are obtained based on the first number 
in $\widehat{E}_k$ and in $\widehat{E}_{\ell}$ using $n_1-1$ doublings.

The elements from $E_{k}$ and $E_{\ell}$ are not inverse to each other, 
and the products of the first element from $E_{k}$ with all elements from $E_{\ell}$ 
are elements with the following degrees
\begin{equation}\label{Produk1}
r_{k,1}+r_{{\ell},1},\; r_{k,1}+2r_{{\ell},1},\cdots,\; r_{k,1}+2^{n_1-1}r_{{\ell},1}\; .
\end{equation}
The elements with the thus obtained $n_1$ degrees belong to certain invariant 
sets $E_j$, $1\le j\le n_1$, which need not necessarily be different.

Next, we consider the products of the second element from $E_{k}$ 
with all elements from $E_{{\ell}}$. But for the elements in $E_{\ell}$  
we want to use the degrees 
$2r_{{\ell},1},4r_{{\ell},1},\cdots, 2^{n_1}r_{{\ell},1}$.
These are also the appropriate degrees for the elements in $E_{\ell}$, 
because they are determined,
based on the preceding second number $2r_{{\ell},1}$, using  $n_1-1$ doublings.
So we obtain  the degrees
\begin{equation}\label{Produk2}
2r_{k,1}+2r_{{\ell},1},\; 2r_{k,1}+4r_{{\ell},1},\cdots,\; 2r_{k,1}+2^{n_1}r_{{\ell},1}\; .
\end{equation}
The values in
(\ref{Produk2})
are twice as large as those in 
(\ref{Produk1}).
With this doubling of the degrees, the corresponding elements remain 
in the same invariant sets $E_j$, $1\le j\le n_1$. 

In the next step, we analogously obtain the degrees
\begin{equation}\label{Produk3}
4r_{k,1}+4r_{{\ell},1},\; 4r_{k,1}+8r_{{\ell},1},\cdots,\; 4r_{k,1}+2^{n_1+1}r_{{\ell},1}\; .
\end{equation}
Compared to
(\ref{Produk2})
we again have a doubling of the degrees, and the corresponding elements belong
to the same invariant sets $E_j$, $1\le j\le n_1$.

This must be done for each of the $n_1$ elements from $E_k$. 
For each element determined by a degree from 
(\ref{Produk1}),
we thus obtain $n_1 - 1$ doublings of degrees and, consequently, all $n_1$ elements 
of the corresponding invariant sets.
The sum of the elements, obtained in the calculation of $S_k\cdot S_{\ell}$, 
is then the sum of the $n_1=2^{\nu+1}$ invariant blocks  $S_j$
to which the elements with the degrees from
(\ref{Produk1})
belong.
$\blacksquare$\newline

The presented proof shows how the product can be calculated practically.
We do not have to multiply every element from the first block 
by every element from the second and check to which invariant set 
the result belongs. We can choose the degree of an element from one of these blocks
and calculate the sums of this degree with the degrees of all elements 
from the other block (calculated {\it modulo }$n$), 
and for thus obtained  $2^{\nu+1}$ numbers check 
to which invariant sets $\widehat{E}_j$ they belong.

The multiplication of invariant blocks is obviously commutative
and we can choose from which of them we will use all elements and from which only one element. 
We can also freely choose which single element we want to use.

The sets $\Bar{E}_j$ are half the size of the sets $\widehat{E}_j$,
and it would be preferable to work with them when determining 
the appropriate invariant sets.
This can be done because an element belongs to the same invariant set as the pair 
to which the element belongs. The degrees of the elements in $E_{\ell}$ are 
$s_{\ell,m}$, $1 \le m \le n_1/2$, and $n-s_{\ell,m}$, $1 \le m \le n_1/2$,
and instead of the degrees in 
(\ref{Produk1}),
we can use the degrees $s_{k,1} + s_{\ell,m}$, $1 \le m \le n_1/2$, and 
$s_{k,1} + n - s_{\ell,m}$, $1 \le m \le n_1/2$.
The elements with these degrees obviously belong to the pairs with the numbers 
$\mathrm{min}(s_{k,1}+s_{\ell,m}, n-(s_{k,1}+s_{\ell,m}))$,\linebreak
 $1\le  m\le n_1/2$, 
and $|s_{k,1}-s_{\ell,m}|$, $1\le  m\le n_1/2$. For these pair numbers, 
we must therefore check to which  $\Bar{E}_j$ they belong.

In the case $n=257$ the invariant sets $E_j$ to which the pairs belong 
can be determined manually, as there are only $16$ sets $\Bar{E}_j$, 
each with $8$ numbers. In the case $n=65537$ it is advisable to use 
computational tools. However, Excel, for example, is sufficient.

In the case $n=257$ we obtain, for example, 
\begin{eqnarray}\label{G1G9}
S_1 S_9=2S_1+2S_3+S_5+2S_6+S_7+2S_9+2S_{11}+S_{13}+2S_{14}+S_{15}\; .
\end{eqnarray}
In the case $n=65537$ we obtain
\begin{eqnarray}\label{G1G1025}
&&S_1\cdot S_{1025} =  2\cdot (S_1+S_{1025})+(S_{24}+S_{1048}) + 2\cdot(S_{155}+S_{1179})\;\;\;\; \\
&&\;\;+(S_{185}+S_{1209}) + (S_{309}+S_{1333}) + (S_{360}+S_{1384}) +(S_{531}+S_{1555})\;\;\;\;  \nonumber \\
&&\;\; + (S_{667}+S_{1691}) + (S_{719}+S_{1743}) + (S_{734}+S_{1758}) \nonumber \\
&&\;\; +2\cdot(S_{778}+S_{1802}) + (S_{841}+S_{1865})+(S_{946}+S_{1970}) \nonumber\; .
\end{eqnarray}

\section{Splitting method and shift property}

Let $n$, $n>17$, be a Fermat prime. Starting with $S$,
\[
S=S_1+S_2+S_3+S_4+\cdots+S_{n_2-1}+S_{n_2}\; ,
\]
we want to successively determine smaller parts so that they consist of invariant blocks.
We will do it in such a way that, starting with the first or second block of the sum 
being splitted, every second block will be consequently added to the new sum.
For $S$ we obtain the parts $F(1,2)$ and $F(2,2)$: 
\begin{eqnarray*}
F(1,2)&=&S_1+S_3+S_5+\cdots+S_{n_2-1}\; ,\\
F(2,2)&=&S_2+S_4+S_6+\cdots+S_{n_2}\; .
\end{eqnarray*}
For $F(1,2)$ and $F(2,2)$ the splitting will be continued. We split $F(1,2)$ 
into the parts $F(1,4)$ and $F(3,4)$ and split $F(2,2)$ into $F(2,4)$ and $F(4,4)$: 
\begin{eqnarray*}
F(1,4)&=&S_1+S_5+S_9+\cdots+S_{n_2-3}\; ,\\
F(3,4)&=&S_3+S_7+S_{11}+\cdots+S_{n_2-1}\; ,\\
F(2,4)&=&S_2+S_6+S_{10}+\cdots+S_{n_2-2}\; ,\\
F(4,4)&=&S_4+S_8+S_{12}+\cdots+S_{n_2}\; .
\end{eqnarray*}
Thus, $S = F(1, 1)$, and for $0\le m\le 2^{\nu}-(\nu+2)$ 
the sets $F(j,2^m)$, $1\le j\le 2^{m}$, are split into  $F(j,2^{m+1})$
and $F(j+2^{m}, 2^{m+1})$,  and the following holds:
\begin{equation}\label{F_Teilung}
F(j,2^m)=F(j,2^{m+1})+F(j+2^m,2^{m+1})\; .
\end{equation}
In the case $m=2^{\nu}-(\nu+2)$ we obtain the individual invariant blocks.

In the notation $F(j,2^m)$ the first parameter $j$, $1\le j\le 2^m$, 
is the number of the invariant block with which the sum begins. The second parameter $2^m$ 
is the amount by which the numbers of invariant blocks increase.\newline

The splitting of the invariant blocks $S_k$ itself we will do analogously.\linebreak
 If $j_1=s_{k,1},\,j_2=s_{k,2},\,\cdots,\,j_{L-1}=s_{k,L-1},\,j_L=s_{k,L}$ 
are the indices of the\linebreak
 $L=n_1/2$ naturally ordered pairs in $E_k$,
we split $S_k=p_{j_1}+p_{j_2}+\cdots+p_{j_{L-1}}+p_{j_L}$ 
into $S_k(1,2)$ and $S_k(2,2)$:
\begin{eqnarray*}
S_k(1,2)&=&p_{j_1}+p_{j_3}+\cdots+p_{j_L-3}+p_{j_L-1}\; ,\\
S_k(2,2)&=&p_{j_2}+p_{j_4}+\cdots+p_{j_L-2}+p_{j_L}\; .
\end{eqnarray*}
Then we split $S_k(1,2)$ into $S_k(1,4)$ and $S_k(3,4)$, and split $S_k(2,2)$ 
into $S_k(2,4)$ and $S_k(4,4)$: \newpage
\begin{eqnarray*}
S_k(1,4)&=&p_{j_1}+p_{j_5}+\cdots+p_{j_L-7}+p_{j_L-3}\; ,\\
S_k(3,4)&=&p_{j_3}+p_{j_7}+\cdots+p_{j_L-5}+p_{j_L-1}\; ,\\
S_k(2,4)&=&p_{j_2}+p_{j_6}+\cdots+p_{j_L-6}+p_{j_L-2}\; ,\\
S_k(4,4)&=&p_{j_4}+p_{j_8}+\cdots+p_{j_L-4}+p_{j_L}\; .
\end{eqnarray*}
Here we have  $S_k=S_k(1,1)$ and for $0\le m\le \nu-1$ we split the sum
$S_k(j,2^m)$, $1\le j\le 2^m$, into $S_k(j,2^{m+1})$ and $S_k(j+2^m,2^{m+1})$:
\begin{equation}\label{S_Teilung}
S_k(j,2^m) = S_k(j,2^{m+1})+S_k(j+2^m,2^{m+1})\; .
\end{equation}
In the case $m=\nu-1$ we obtain the individual pairs.

The parameters of $S_k(j,2^m)$ are as follows: $k$ is the block number,\linebreak
$j$ indicates that the sum begins with the $j$-th pair in the ordered set $E_k$, 
and $2^m$ is the step with which the pairs in $E_k$ are selected.

In the case $n=17$, the splittings were carried out in precisely this way.
We could see that it is the only suitable variant in this case.\newline

{\bf Remark. }The presented splittings are analogous to the splittings for Gaussian periods.
However, the author came to the splitting approach in a completely different way
on the basis of invariant sets.\newline

The ordered invariant sets are formed by multiplying the degrees of elements 
by $2^{h_1}$ and $g^{h_2}$, where $g$ is a primitive root and $h_1$ and $h_2$
are non-negative integers.
We want to show what happens to the parts $F(j, 2^m)$ and $S_k(j, 2^m)$ 
when the degrees of the elements are multiplied by $2^{h_1} g^{h_2}$.
To state the result more simply, we will use the following function 
for positive integers $k_1$ and $k_2$: 
\begin{equation}\label{rho}
\rho(k_1,k_2)=\begin{cases}
k_2, & \texttt{if }\mathrm{rest}(k_1,k_2)=0\\
\mathrm{rest}(k_1,k_2), & \texttt{if }\mathrm{rest}(k_1,k_2)>0
\end{cases}
\end{equation}

If we multiply the degrees of the elements by $2^{h_1}$, we obtain a cyclic shift 
for the corresponding $n_1$ elements in each invariant set $E_k$.
The sum $S_k$ does not change in this process, and this means
that multiplication by $2^{h_1}$ plays no role for $F(j, 2^m)$.
Multiplying the degrees by $g^{h_2}$ leads to the cyclic shift for the $n_2$ 
invariant sets $E_1, \dots, E_{n_2}$ of the height $h_2$. 
The first summand $S_j$ in $F(j, 2^m)$ goes then to $S_{\rho(j+h_2, n_2)}$.
This shift does not change the distance $2^m$ between the numbers 
of the invariant sets and,
consequently, $F(j,2^m)$ transforms into a $F(j_1,2^m)$, where $1\le j_1\le 2^m$. 
Since the distance between $j_1$ and $j + h_2$ is a multiple of $2^m$, 
it is obvious that  $j_1 = \rho(j + h_2, 2^m)$. Thus, $F(j, 2^m)$ transforms 
into $F(\rho(j+h_2,2^m),2^m)$.

If we multiply the degrees of the elements by $2^{h_1}$, we obtain a cyclic shift
with the height $h_1$ for the $n_1/2$ numbers $s_{k,1}, \cdots, s_{k,n_1/2}$ of the pairs 
from each invariant set $E_k$. 
The summands in $S_k(j, 2^m)$ have the numbers 
$s_{k,j}, s_{k,j+2^m}, s_{k,2\cdot 2^m},\cdots$, 
and $2^m$ is the amount by which the index $j$ increases.
For the cyclic shift this amount remains unchanged and we obtain the values
$s_{k,j_1}, s_{k,j_1+2^m}, s_{j_1+2^m},\cdots$, where $1 \le j_1 \le 2^m$.
Analogously to the previous we then obtain $j_1 = \rho(j + h_1 , 2^m )$,
and $S_k (j, 2^m )$ transforms into\linebreak
$S_k (\rho(j + h_1 , 2^m ), 2^m )$.
If we  multiply the degrees
by $g^{h_2}$, the set $E_k$ becomes $E_{\rho(k+h_2 ,n2)}$ and the order
of the pairs does not change.
This means that the multiplication of the degrees by $2^{h_1} g^{h_2}$, 
transforms the value $S_k(j, 2^m)$ into $S_{\rho(k+h_2, n_2)}(\rho(j+h_1, 2^m), 2^m)$.

To split $S$, we need the products of the corresponding parts.
We want to present a feature of these products, which we will denote shift property 
in the following.

We know the rule governing how the parameters of the parts $F(j, 2^m)$ and $S_k(j, 2^m)$ 
change when the degrees of the elements are multiplied by $2^{h_1} g^{h_2}$.
If we have a change from the parameters of $U$ 
to the parameters of $V$ just as it must be due to multiplication of degrees 
by $2^{h_1} g^{h_2}$, we say that we come from $U$ to $V$
via shift with the height $\textit{shift}_1 = h_1$ for the pair numbers
and $\textit{shift}_2 = h_2$ for the block numbers.
\newline

{\bf Proposition 4. }Let $X_1$ and $X_2$ be certain parts from the proposed 
split of $S$ and let the product $X_1 X_2$ be representable as 
$X_1 X_2 = P(U_1,\cdots,  U_m)$
with a polynomial $P(x_1,\cdots,x_m)$ with real coefficients 
and proposed parts $U_1,\cdots,U_m$ from $S$. 
If we come from $X_1$ and $X_2$ to $Y_1$ and $Y_2$ via the same shift, 
then the product $Y_1 Y_2$ can be represented as $Y_1 Y_2=P(V_1,\cdots,V_m)$,
where $V_1,\cdots,V_m$ are such 
that we come from each $U_k$ to $V_k$, $1\le k\le m$, via the same shift 
as that from $X_1$ and $X_2$ to $Y_1$ and $Y_2$.

{\bf Proof. }Let $\gamma=2^{h_1}g^{h_2}$ be the factor by which the degrees 
are multiplied during the shift.
$X_1$ and $X_2$ can be represented as follows:
$X_1=\sum_{k=1}^{n-1}a_{1,k}z^k$ and $X_2=\sum_{k=1}^{n-1}a_{2,k}z^k$,
and for the product $X_1 X_2$, we then obtain the representation
$X_1 X_2=c_0+\sum_{k=1}^{n-1}c_kz^k$.
Here, the coefficients $a_{1,k}$ and $a_{2,k}$ are $0$ or $1$, and $c_k$,
$0 \le k \le n-1$, are non-negative integers.
The constant $c_0$ we obtain here using the products of inverse elements 
in $X_1$ and $X_2$.

For $Y_1$ and $Y_2$ the following representations then apply:
$Y_1=\sum_{k=1}^{n-1}a_{1,k}z^{\gamma\cdot k}$, 
$Y_2=\sum_{k=1}^{n-1}a_{2,k}z^{\gamma\cdot k}$,
and it is evident that
$Y_1 Y_2=c_0+\sum_{k=1}^{n-1}c_kz^{\gamma\cdot k}$.

But if $X_1 X_2=P\big{(}\sum_{k=1}^{n-1}b_{1,k}z^k,\cdots,
\sum_{k=1}^{n-1}b_{m,k}z^k\big{)}=c_0+\sum_{k=1}^{n-1}c_kz^k$
holds, then we obviously obtain:
\[
\it{ }\;\;\;\;\;\;\; 
P\Big{(}\sum_{k=1}^{n-1}b_{1,k}z^{\gamma\cdot k},\cdots,
\sum_{k=1}^{n-1}b_{m,k}z^{\gamma\cdot k}\Big{)}= 
c_0+\sum_{k=1}^{n-1}c_kz^{\gamma\cdot k}=Y_1  Y_2\;.\;\;
\;\;\blacksquare
\]
\section{Feasibility of the splitting steps}

We want to show that the $2^{\nu}-1$ splitting steps used to come from $S$
to the individual pairs can be realized with compass and straightedge.
The following two propositions show that, for any Fermat prime $n$, 
the regular $n$-gon can be constructed.
The proofs also show how the construction can be practically realized.\newline

{\bf Proposition 5. }The product $F(j,2^{m+1})\cdot F(j+2^m, 2^{m+1})$
is a linear combination of the values $F(k,2^m)$, $1\le k\le 2^m$, with non-negative
integer coefficients and can be constructed on the basis of these values. 

{\bf Proof. }We consider at first the product $2F(1,2^{m+1})\cdot F(1+2^m, 2^{m+1})$.
$F(1,2^m)$ consists of $q_m=n_2/2^m$ summands $S_{1+k\cdot 2^m}$, $0\le k\le q_m-1$,
and for $F(1,2^{m+1})$ and $F(1+2^m,2^{m+1})$ we have: 
\begin{eqnarray*}
&F(1,2^{m+1})=S_1+S_{1+2\cdot 2^{m}}+S_{1+4\cdot 2^{m}}+S_{1+6\cdot 2^{m}}
+\cdots+S_{1+(q_m-2)2^{m}}\; ,\\
&\!\!\!\!F(1+2^m,2^{m+1})=S_{1+2^m}+S_{1+3\cdot 2^{m}}+S_{1+5\cdot 2^{m}}+S_{1+7\cdot 2^{m}}
+\cdots+S_{1+(q_m-1)2^{m}}\; .\;\;
\end{eqnarray*}

We will use the cyclic property of the summands in $F(1,2^{m+1})$ and in $F(1+2^m,2^{m+1})$
and calculate $2F(1,2^{m+1})\cdot F(1+2^m, 2^{m+1})$ by alternately 
taking the summands  $S_{1+k\cdot 2^m}$ 
from $F(1,2^{m+1})$ and $F(1+2^m,2^{m+1})$ and multiplying them with all summands
of the part, either $F(1,2^{m+1})$ or $F(1+2^m,2^{m+1})$, 
to which $S_{1+k\cdot 2^m}$ does not belong.
In this part the ordering of the summands should be chosen so that it begins with the term
$S_{1+(k+1)\cdot 2^m}$. So we obtain:
\begin{eqnarray*}
&&2F(1,2^{m+1})\cdot F(1+2^m, 2^{m+1})=\\
&&\;\;S_{1+0\cdot 2^m} S_{1+1\cdot 2^m}\;+S_{1+0\cdot 2^m} S_{1+3\cdot 2^{m}}+\cdots 
+ S_{1+0\cdot 2^m} S_{1+(q_m-1)2^{m}}\nonumber\\
&&+S_{1+1\cdot 2^m} S_{1+2\cdot 2^m}+S_{1+1\cdot 2^m} S_{1+4\cdot 2^{m}}
+\cdots+ S_{1+1\cdot 2^m} S_{1+q_m 2^{m}}\nonumber\\
&&+S_{1+ 2\cdot 2^m} S_{1+3\cdot 2^m}+S_{1+ 2\cdot 2^m} S_{1+5\cdot 2^{m}}
+\cdots+ S_{1+ 2\cdot 2^m} S_{1+(q_m+1) 2^{m}}\nonumber\\
&&+S_{1+3\cdot 2^m} S_{1+4\cdot 2^m}+S_{1+3\cdot 2^m} S_{1+6\cdot 2^{m}}+
\cdots+ S_{1+3\cdot 2^m} S_{1+(q_m+2) 2^{m}} \nonumber\\
&&+\cdots\nonumber
\end{eqnarray*}
Here, we show only $4$ of the total $q_m$ rows in detail.

In this representation, the number of the first factor remains constant 
within the rows and always increases by the amount $2^m$ when we go to the next row.
The numbers of the second factors differ within the rows, 
but increase by the same amount  $2^m$ when we go to the next row.
Due to the shift property, the numbers of the blocks obtained with the products 
also increase by the amount $2^m$ when we go to the next row.
With the $q_m/2$ products in the first row
\begin{eqnarray}\label{F_Zeile_1}
S_1 S_{1+2^m}+S_1 S_{1+3\cdot 2^{m}}+\cdots+S_1 S_{1+(q_m-3)2^{m}}+S_1 S_{1+(q_m-1)2^{m}}
\end{eqnarray}
we obtain invariant blocks, and these blocks belong to certain values 
$F(k,2^m)$, $1\le k\le 2^m$.
But as the block numbers increase in the following $q_m-1$ rows repeatedly
by the amount  $2^m$, we obtain all summands from these values $F(k,2^m)$.
This means that the product $2F(1,2^{m+1})\cdot F(1+2^m, 2^{m+1})$
is a sum of values $F(k,2^m)$, and these are the values to which the blocks 
obtained with the products in the first row 
(\ref{F_Zeile_1})
belong.

We want to show that using the first half of the products in 
(\ref{F_Zeile_1}),
we obtain the same values $F(k,2^m)$ as using the second half of the products.
To do this, we compare 
the products $S_1S_{1+(2k-1)\cdot 2^m}$ and 
$S_1S_{1+(q_m-(2k-1))2^m}$, $1\le k\le q_m/4$.
We consider a shift with the height $\textit{shift}_2=(q_m-(2k-1))2^m$\linebreak  
for the block numbers.
Via this shift  $S_1$ goes to $S_{1+(q_m-(2k-1))2^m}$ and\linebreak
 $S_{1+(2k-1)2^m}$ 
goes to $S_{1+q_m\cdot 2^m}=S_{1+(n_2/2^m)\cdot 2^m}=S_{1+n_2}=S_1$.
Hence this shift turns $S_1S_{1+(2k-1)\cdot 2^m}$ into $S_1S_{1+(q_m-(2k-1))2^m}$. 
However, since the shift height is a multiple of $2^m$, the blocks in the product
$S_1S_{1+(2k-1)\cdot 2^m}$ belong to the same values $F(k,2^m)$
as the blocks in product $S_1S_{1+(q_m-(2k-1))2^m}$ .

We thus see that $F(1,2^m)\cdot F(1+2^m)$ is a linear combination of the values $F(k,2^m)$
with integer coefficients $\mu_{k,2^m}$ 
\begin{eqnarray}\label{F1Fm_aus_0,5Zeile}
F(1,2^{m+1})\cdot F(1+2^m, 2^{m+1})=\sum_{1\le k\le 2^{m}} \mu_{k,2^m} F(k,2^m)\; ,
\end{eqnarray}
and it can, in fact, be determined using the first half of the products\linebreak
in 
(\ref{F_Zeile_1}).
The coefficient $\mu_{k,2^m}$ can be easily calculated here because all the blocks $S_v$ 
determined by the products in
(\ref{F_Zeile_1})
can be calculated, and each individual block $S_v$ belongs to $F(k,2^m)$ with $k=\rho(v,2^m)$.

We come from $F(1,2^{m+1})$ and $F(1+2^m, 2^{m+1})$ to $F(j,2^{m+1})$ and\linebreak
 $F(j+2^m, 2^{m+1})$
via shift with the height $\textit{shift}_2=j-1$ for the block numbers. 
Due to the shift property we then obtain: 
\begin{eqnarray}\label{FjFj+2^m}
F(j,2^{m+1})\cdot F(j+2^m, 2^{m+1})=\sum_{1\le k\le 2^m} \mu_{k,2^m}
F(\rho(k+j-1,2^m),2^m)\; .
\end{eqnarray}
This linear combination certainly 
can be constructed with a compass and straightedge if the required values
$F(k,2^m)$ are known.
$\blacksquare$

{\bf Remark. }The product $F(1,2)\cdot F(2,2)=\mu_{1,1}\cdot S$ 
represents an even covering of $S$.
$F(1,2)$ and $F(2,2)$ have $(n-1)/2$ elements, $F(1,2)\cdot F(2,2)$ 
obviously has $(n-1)^2/4$ elements, and the height of the covering 
for the product $F(1,2)\cdot F(2,2)$ is then $(n-1)/4$.
Therefore:
\begin{equation}\label{F1F2_einfach}
F(1,2)\cdot F(2,2)=(n-1)/4\cdot S=-(n-1)/4\; .
\end{equation}

In the case $n=257$ we can easily show that $F(1,4)\cdot F(3,4)$ is also a constant. 
Indeed, in this case we can analytically calculate $S_1S_3+S_1S_7$ and
determine the amount of invariant blocks from $F(1,2)$ and $F(2,2)$ in this sum.
We thus obtain $\mu_{1,2}=\mu_{2,2}=16=(n-1)/16$, and therefore
\begin{equation}\label{F14F34}
F(1,4)\cdot F(3,4)=16\cdot (F(1,2)+F(2,2))= 16\cdot S =-16=-(n-1)/16\; .
\end{equation}

In the case $n = 65537$, the product $F(1, 4)\cdot F(3, 4)$ cannot be calculated 
analytically so easily, because it results in $256$ products of blocks $S_1\cdot S_j$.
However, the numerical calculation yields $\mu_{1,2}=\mu_{2,2}=4096=(n-1)/16$, 
and therefore $F(1,4)\cdot F(3,4)=-(n-1)/16$ also holds true  for $n = 65537$.

In both cases, i.e. for $n=257$ and  $n=65537$,
we come from  $F(1,4)$  and $F(3,4)$ to $F(2,4)$ and $F(4,4)$ 
via shift with the height $\textit{shift}_2=1$ 
for the block numbers, and therefore we also have $F(2,4)\cdot F(4,4)=-(n-1)/16$.\newline

{\bf Proposition 6.} The product $S_k(j,2^{m+1})\cdot S_k(j+2^m,2^{m+1})$
can be calculated and constructed with a compass and straightedge
on the basis of the values $S_{\ell}(q,2^m)$, $1\le {\ell}\le n_2$, $1\le q\le 2^m$.

{\bf Proof. }$S_{\ell}(q,2^m)$ has $L_m=2^{\nu}/2^m$ pairs as summands, 
and the numbers of the pairs in $S_{\ell}(q,2^m)$, starting 
with the first number $s_{{\ell},q}$, 
increase upon multiplication by the factor $M_m=2^{2^m}$. Therefore we simply have 
\begin{equation}\label{barE_neu}
S_{\ell}(q,2^m)=p_{s_{{\ell},q}}+p_{s_{{\ell},q} M_m}+p_{s_{{\ell},q} M_m^2}
+p_{s_{{\ell},q} M_m^3}
+\cdots +p_{s_{{\ell},q} M_m^{L_m-1}}\; . 
\end{equation}
The calculation {\it modulo }$n$ and the conversion from element number 
to pair number are omitted here.

If the numbers of the pairs in 
(\ref{barE_neu})
are multiplied by the factor $M_m^r$, we obtain the same value $S_{\ell}(q,2^m)$.
Indeed, as a result of the multiplication with the factor $M_m$, the numbers 
$s_{{\ell},q}, s_{{\ell},q} M_m, s_{{\ell},q} M_m^2, s_{{\ell},q} M_m^3,
\cdots,s_{{\ell},q} M_m^{L_m-1}$\linebreak
become the numbers 
$s_{{\ell},q} M_m, s_{{\ell},q} M_m^2, s_{{\ell},q} M_m^3,\cdots,s_{{\ell},q} M_m^{L_m-1}, 
s_{{\ell},q}M_m^{L_m}$\linebreak
and we have $p_{s_{{\ell},q} M_m^{L_m}}=p_{s_{{\ell},q}}$ because 
\[
s_{{\ell},q} M_m^{L_m}=s_{{\ell},q}\cdot(2^{2^m})^{2^{\nu}/2^m}=s_{{\ell},q}\cdot 2^{2^{\nu}}=
s_{{\ell},q}\cdot(n-1)\equiv -s_{{\ell},q} \pmod{n}\; .
\]
Multiplication by $M_m$ thus leads to a cyclic shift for the pairs in 
(\ref{barE_neu}) :
the old second pair becomes the new first pair.
Multiplication by $M_m^r$ results in $r$ cyclic shifts.
The sum $S_{\ell}(q,2^m)$ does not change in this process.

We can see, among other things, that when calculating $S_{\ell}(q,2^m)$, 
we can start with any of its pairs.
However,  $L_m-1$ additional pairs must subsequently be included as summands,
and the numbers of these pairs must be determined successively by multiplication 
with the factor  $M_m$.

It is also clear that any sum of $L_m$ pairs, in which the numbers 
are calculated successively by multiplication with the factor $M_m$,
is a sum $S_{\ell}(q,2^m)$.
The parameters ${\ell}$ and $q$ of the sum $S_{\ell}(q,2^m)$ formed in this way 
can be determined as follows. For the number 
$k_0=\min(\mathrm{rest}(s,n),n-\mathrm{rest}(s,n))$
of the first pair, $1\le k_0\le (n-1)/2$, 
we determine to which set $\Bar{E}_{\ell}$ it belongs.
If $k_0$ is the $v$-th number in $\Bar{E}_{\ell}$, i.e. $k_0=s_{{\ell},v}$,
then we obviously have the sum $S_{\ell}(q,2^m)$ with  $q=\rho(v,2^m)$.

We consider first the product $S_1(1,2^{m+1})\cdot S_1(1+2^m,2^{m+1})$.\linebreak
$S_1(1,2^{m+1})$ and $S_1(1+2^m,2^{m+1})$ consist of $L_{m+1}=L_m/2$ pairs,
and the numbers of these pairs grow by the factor $M_{m+1}=2^{2^{m+1}}=M_m^2$.
Therefore $S_k(1,2^{m+1})$ and $S_k(1+2^m,2^{m+1})$  can be represented as follows:
\begin{eqnarray*}
S_1(1,2^{m+1})&=&p_1+p_{M_m^2}+p_{M_m^4}+p_{M_m^6}+\cdots+p_{M_m^{L_m-2}}\; ,\\
S_1(2^m,2^{m+1})&=&p_{M_m}+p_{M_m^3}+p_{M_m^5}+p_{M_m^7}+\cdots+p_{M_m^{L_m-1}}\; .
\end{eqnarray*}

We want to calculate the product $S(1,2^{m+1})\cdot S(1+2^m,2^{m+1})$ using the squares.
Because of
(\ref{S_Teilung})
we have:
\begin{eqnarray}\label{2S1S2}
&&2S_1(1,2^{m+1})\cdot S_1(1+2^m,2^{m+1})=(S_1(1,2^{m+1})+ 
S_1(1+2^m,2^{m+1}))^2\nonumber\\
&&-S_1^2(1,2^{m+1})-S_1^2(1+2^m,2^{m+1})=S_1^2(1,2^m)-Q_m-P_m\; ,
\end{eqnarray}
where  $Q_m$ is the sum of the squares of all elements from $S_1(1,2^{m+1})$
and from $S_1(1+2^m,2^{m+1})$, and $P_m$ is twice the sum of the products 
of distinct pairs from $S_1(1,2^{m+1})$ and from $S_1(1+2^m,2^{m+1})$.

$Q_m$ is the sum of the squares of the $L_m$ pairs from $S_1(1,2^{m})$, 
and the following holds: 
\begin{eqnarray}\label{Q}
Q_m&=&p^2_1+p^2_{M_m}+p^2_{M_m^2}+p^2_{M_m^3}+\cdots+p^2_{M_m^{L_m-1}}\nonumber\\
&=&\left(p_2+p_{2M_m}+p_{2M_m^2}+p_{2M_m^3}+\cdots+p_{2M^{(L_m-1)}}\right)
+2L_m\nonumber\\
&=& S_1(2,2^m)+2^{\nu+1-m}\; .
\end{eqnarray}

At the final split, i.e. at $m=\nu-1$, we have
only one pair in $S_1(1,2^{m+1})$ and in $S_1(1+2^m,2^{m+1})$, and therefore $P_{\nu-1}=0$.
In case $m<\nu-1$ we can calculate the double sum of products 
of distinct pairs from $S_1(1,2^{m+1})$ and from $S_1(1+2^m,2^{m+1})$
by alternately taking the pairs $p_{M^k}$ from $S_1(1,2^{m+1})$ and $S_1(1+2^m,2^{m+1})$
and multiplying them by the other pairs from the part, 
either $S_1(1,2^{m+1})$ or $S_1(1+2^m,2^{m+1})$, to which $p_{M^k}$ belongs.
In this part the ordering of the summands should be chosen so that
$p_{M^k}$ is the first pair.
In this way we obtain the following representation for $P_m\,$:
\begin{eqnarray}\label{P_m}
P_m&=&\;(\;\; p_1p_{M_m^2}\;\;+\;\;p_1p_{M_m^4}\;\;+\;p_1p_{M_m^6}\;+\cdots \,+p_1p_{M_m^{L_m}})\nonumber\\ 
&&+(p_{M_m}p_{M_m^3}+p_{M_m}p_{M_m^5}+p_{M_m}p_{M_m^7}+\cdots +p_{M_m}p_{M_m^{L_m+1}} )\nonumber\\
&&+(p_{M_m^2}p_{M_m^4}+p_{M_m^2}p_{M_m^6}+p_{M_m^2}p_{M_m^8}+\cdots +p_{M_m^2}p_{M_m^{L_m+2}})\\
&&+(p_{M_m^3}p_{M_m^5}+p_{M_m^3}p_{M_m^7}+p_{M_m^3}p_{M_m^9}+\cdots +p_{M_m^3}p_{M_m^{L_m+3}})\nonumber \\
&&+\cdots \nonumber
\end{eqnarray}
In this representation, there are $L_m$ rows --- one for each pair in $S_1(1,2^m)$.
We are showing only $4$ of them in detail.

The number of the first pair is constant within the products of each row in
(\ref{P_m})
and increases by the factor $M_m$ when we go to the next row.
The numbers of the second pairs also increase by the factor $M_m$ 
when we go to each next row.
But if the numbers of both pairs increase by the factor $M_m$, then the numbers 
of the two pairs determined by the product also increase by the factor $M_m$.
Thus, we see that for the  $L_m-2$ pairs determined by the $L_m/2-1$ products 
in the first row
\begin{eqnarray}\label{Zeile1}
p_1p_{M_m^2}+p_1p_{M_m^4}+p_1p_{M_m^6}+\cdots+p_1p_{M_m^{L_m/2}} +
p_1p_{M_m^{L_m/2+2}}+\cdots +p_1p_{M_m^{L_m}}\;\;
\end{eqnarray}
the numbers increase by the factor $M_m$ in the subsequent $L_m-1$ rows,
and these summands represent all pairs of certain values $S_{\ell}(q,2^m)$.
Therefore, $P_m$ is the sum of $L_m-2$ values $S_{\ell}(q,2^m)$.

Determining the values $S_{\ell}(q,2^m)$ to which the pairs belong is not entirely simple,
and we will examine more closely the pairs obtained using
(\ref{Zeile1})
We consider the case $1\le m<\nu-1$ and have in 
(\ref{Zeile1})
an odd number of products. We will compare the products that are located 
symmetrically about the center.
These are $p_1p_{M_m^{2k}}$ and $p_1p_{M_m^{L_m-2k}}$, and we obtain the pairs
with the numbers
$M_m^{2k}-1$ and $M_m^{2k}+1$, as well as  $M_m^{L_m-2k}-1$ and $M_m^{L_m-2k}+1$.
Here we have $M_m=2^{2^m}$ and $L_m=2^{\nu}/2^m$, and 
if we multiply the index of a pair from $S_{\ell}(q,2^m)$ by the factor $M_m^{2k}$,
we obtain a pair in the same $S_{\ell}(q,2^m)$.
It holds
\[
M_m^{L_m}=(2^{2^m})^{2^{\nu}/2^m}=2^{2^{\nu}}=n-1\equiv -1\pmod{n}\; ,
\]
and therefore
\begin{eqnarray*}
&&M_m^{2k}\cdot(M_m^{L_m-2k}-1)=-M_m^{2k}+M_m^{L_m}\equiv -(M_m^{2k}+1)\pmod{n}\; ,\\
&&M_m^{2k}\cdot(M_m^{L_m-2k}+1)=M_m^{2k}+M_m^{L_m}\equiv M_m^{2k} -1\pmod{n}\; .
\end{eqnarray*}
The symmetrically positioned products thus yield pairs from the same $S_{\ell}(q,2^m)$.

In the middle, we have the product $p_1p_{M_m^{L_m/2}}$, and we obtain the pairs 
with the numbers $k_1=2^{2^{\nu-1}}-1$ and $k_2=2^{2^{\nu-1}}+1$.
We want to show that the pairs $p_{k_1}$ and $p_{k_2}$
belong to the same $S_{\ell}(q,2^m)$.
Indeed, it holds:
\[
2^{2^{\nu-1}}\cdot(2^{2^{\nu-1}}+1)=2^{2^{\nu-1}}+2^{2^{\nu}}=2^{2^{\nu-1}}+(n-1)
\equiv 2^{2^{\nu-1}}-1\pmod{n}
\]
and thus we come from the pair $p_{k_2}$ to $p_{k_1}$ by $2^{\nu-1}$ 
doublings of the number $k_2$.
This means that $p_{k_1}$ and $p_{k_2}$ belong to the same invariant block $S_{\ell}$.
But for $m < \nu - 1$ the following still holds:
\[
2^{2^{\nu-1}}=(2^{2^m})^{2^{\nu-1}/2^m}=M_m^{2^{\nu-1-m}}
\]
and we come from $k_2$ to $k_1$ via multiplying $k_2$ by the factor $M_m$
a total of $2^{\nu-1-m}$ times.
This means that $p_{k_1}$ and $p_{k_2}$ belong to the same $S_{\ell}(q,2^m)$.

We thus see that to calculate $P_m$ for $m < \nu - 1$ we indeed have to
determine the appropriate values $S_{\ell}(q, 2^m)$ only for the products of pairs 
positioned left from the middle summand $p_1p_{M_m^{L_m/2}}$ in
(\ref{Zeile1})
and for one pair from the product in this middle summand.
We denote the sum of these values by $H_m$ in the following.
In this notation, $H$ is intended to signify "half". Thus, $P_m = 2H_m$, and due to
(\ref{2S1S2})
it holds
\begin{eqnarray}\label{S1S2}
S_1(1,2^{m+1})\cdot S_1(2^m,2^{m+1})=\frac{S_1^2(1,2^m)-S_1(2,2^m)}{2}-2^{\nu-m}+H_m\; .
\end{eqnarray}

We come from $S_1(1,2^{m+1})$ and $S_1(1+2^m,2^{m+1})$ to $S_k(j,2^{m+1})$ and\linebreak
$S_k(j+2^m,2^{m+1})$ via  shift with the height $\textit{shift}_1=j-1$ for the pair numbers 
and $\textit{shift}_2=k-1$ for the block numbers.
Due to the shift property we obtain then $S_k(j,2^{m+1})\cdot S_k(j+2^m,2^{m+1})$
on the basis of
(\ref{S1S2})
using the same shift.
The product $S_k(j,2^{m+1})\cdot S_k(j+2^m,2^{m+1})$ can also be constructed
if all required values $S_{\ell}(q,2^m)$ are known.
$\blacksquare$ 

{\bf Remark. }We could also simply take the pairs from $S_k(j, 2^{m+1})$ one by one 
and multiply them by all pairs in $S_k(j+2^m, 2^{m+1})$.
In doing so, we could repeatedly choose the ordering for $S_k(j+2^m,2^{m+1})$
in such a way that, analogously to 
(\ref{P_m}),
we obtain sums of products that consist of values $S_{\ell}(q,2^m)$.
However, this approach requires more values of $S_{\ell}(q,2^m)$ 
than the approach using squares.

\section{Regular 257-gon}

{\bf Step 1. }We split $S$ into $F(1,2)$ and $F(2,2)$.
Here the sum and the product are known: $F(1,2)+F(2,2)=-1$ and 
$F(1,2)\cdot F(2,2)=-(n-1)/4=-64$. Therefore $F(1,2)$ and $F(2,2)$ are 
the solutions of the the quadratic equation 
\begin{eqnarray}\label{QuadGlStep1}
x^2+x-64=0
\end{eqnarray} 
and due to $F(1,2)>F(2,2)$ we have
\begin{eqnarray*}
F(1,2)=\frac{-1+\sqrt{257}}{2}\; ,\;\;
F(2,2)=\frac{-1-\sqrt{257}}{2}\; .
\end{eqnarray*}

{\bf Step 2.} We split $F(j,2)$ into $F(j,4)$ and $F(2+j,4)$, $1\le j\le 2$.\newline
The sum and the product are known here: $F(j,4)+F(2+j,4)=F(j,2)$, 
$F(j,4)\cdot F(2+j,4)=-(n-1)/16=-16$.
Therefore we can calculate $F(j,4)$ and $F(2+j,4)$ on the basis of the quadratic equation 
\[
x^2-F(j,2)x-16=0\; .
\]
Due to $F(1,4)>F(3,4)$ and $F(2,4)>F(4,4)$ we obtain for $1\le j\le 2$
\begin{eqnarray*}
F(j,4)&=&\frac{F(j,2)+\sqrt{F^2(j,2)+64}}{2}\; ,\\
F(2+j,4)&=&\frac{F(j,2)-\sqrt{F^2(j,2)+64}}{2}\; .
\end{eqnarray*}
It is possible to rearrange these formulas so that $F^2(j,2)$ disappears.
Indeed, as shown in step 1,
(\ref{QuadGlStep1}),
we have for $F(j,2)$ the equality
\[
F^2(j,2)+F(j,2)-64=0
\]
and hence it follows for $1\le j\le 2$
\begin{eqnarray*}
F(j,4)&=&\frac{F(j,2)+\sqrt{128-F(j,2)}}{2}\; ,\\
F(2+j,4)&=&\frac{F(j,2)-\sqrt{128-F(j,2)}}{2}\; .
\end{eqnarray*}

{\bf Step 3.} We consider at first the split of
$ F(1,4)=S_1+S_5+S_{9}+S_{13}$
\newline
into $\;F(1,8)=S_1+S_{9}\;$ and $\;F(5,8)=S_5+S_{13}$.

To calculate the product $F(1,8)\cdot F(5,8)$
we need the the integer coefficients $\mu_{k,4}$ for $F(k,4)$, $1\le k\le 4$.
We obtain here
$\mu_{1,4}=2$, $\mu_{2,4}=5$, $\mu_{3,4}=4$ and $\mu_{4,4}=5$.
Thus we have
\begin{eqnarray*}
F(1,8)\cdot F(5,8)&=&2F(1,4)+5F(2,4)+4F(3,4)+5F(4,4)\\
&=&-5-3F(1,4)-F(3,4)
\end{eqnarray*}
and it follows for $F(j,8)$ and $F(4+j,8)$, $1\le j\le 4\;$:
\begin{eqnarray*}
F(j,8)+F(4+j,8)&=&F(j,4)\\
F(j,8)\cdot F(4+j,8)&=&-5-3F(j,4)-F(\rho(2+j,4),4)
\end{eqnarray*}
The values $F(j,8)$ and $F(4+j,8)$, $1\le j\le 4$, definitely 
can then be calculated and constructed on the basis of the relevant quadratic equations 
provided that
\begin{eqnarray*}
F(1,8)>F(5,8),\; F(2,8)<F(6,8),\; F(3,8)>F(7,8),\; F(4,8)<F(8,8).
\end{eqnarray*}
We do not apply the possible transformations of the solutions so that
the square $F^2(j,4)$ disappears.
This is not absolutely necessary as the square can be constructed
with a compass and straightedge. 

{\bf Step 4.} We consider at first the split of $F(1,8)=S_1+S_{9}$ 
into
\newline
$F(1,16)=S_1$ and $F(9,16)=S_9$.

The product $S_1\cdot S_9$ is already calculated,
(\ref{G1G9}).
A simple transformation of this result yields:
\begin{eqnarray*}
S_1\cdot S_9&=&2F(1,8)+2F(3,8)+F(5,8)+2F(6,8)+F(7,8)\\
&=&-2-F(5,8)-F(7,8)\; .
\end{eqnarray*}

For $S_{j}$ and $S_{8+j}$, $1\le j\le 8$, we have therefore:
\begin{eqnarray*}
S_j+S_{8+j}&=&F(j,8)\\
S_j\cdot S_{8+j}&=&-2-F(\rho(4+j,8),8)-F(\rho(6+j,8),8)
\end{eqnarray*}
For the solutions of the relevant quadratic equations it must be kept in mind that 
\begin{eqnarray*}
S_1>S_8, S_2>S_{10},\; S_3>S_{11},\; S_4>S_{12},\; S_5>S_{13},\; S_6>S_{14}, S_7<S_{15}\; .
\end{eqnarray*}

From here on we consider steps 5–7 in reverse order.
Due to the special interest in the value $p_1$, we consider the steps so 
that we come to $p_1$ without unnecessary splittings.\newline

{\bf Step 7.} The last split to obtain $p_1$ is the split of
 $S_1(1,4)=p_1+p_{16}$ into $p_1$ and $p_{16}$. 
In this case we have no products in
(\ref{Zeile1}),
and therefore $H_2=0$. Hence we obtain:
\begin{eqnarray*}
p_1+p_{16}&=&S_1(1,4)\\
p_1\cdot p_{16}&=&\frac{S^2_1(1,4)-S_1(2,4)}{2}-2
\end{eqnarray*}
The values $p_1$ and $p_{16}$ can be uniquely determined
on the basis of the relevant quadratic equation by taking into account that $p_1>p_{16}$. 

For the calculations in step 7 we need in any case the values
$S_1(1,4)$ and $S_1(2,4)$, and it is necessary to get them
in step 6 at the splitting of $S_1(1,2)$ and $S_1(2,2)$.\newline

{\bf Step 6.} We split at first $S_{1}(1,2)=p_1+p_4+p_{16}+p_{64}$
into $S_1(1,4)=p_1+p_{16}$ and $S_1(3,4)=p_4+p_{64}$.

In this case
(\ref{Zeile1})
consists only of the middle product $p_1p_{16}$.
The pair $p_{15}$ belongs to $S_9(2,2)$ and hence $H_1=S_9(2,2)$.

For $S_1(1,4)$ and $S_1(3,4)$  we have therefore:
\begin{eqnarray*}
S_1(1,4)+S_1(3,4)&=&S_1(1,2)\\
S_1(1,4)\cdot S_1(3,4)&=&\frac{S^2_1(1,2)-S_1(2,2)}{2}-4-S_9(2,2)
\end{eqnarray*}
and in the calculation of $S_1(1,4)$ and $S_1(3,4)$ 
on the basis of the relevant quadratic equation it is to consider that $S_1(1,4)>S_1(3,4)$.

To obtain $S_1(2,4)$, we split $S_1(2,2)$ into $S_1(2,4)$ and $S_1(4,4)$.
Compared to $S_1(1,4)$ and $S_1(3,4)$ we have here a shift 
with the height $\textit{shift}_1=1$ for the pair numbers.
For $S_1(2,4)$ and $S_1(4,4)$  we get therefore:
\begin{eqnarray*}
S_1(2,4)+S_1(2,4)&=&S_1(2,2)\\
S_1(2,4)\cdot S_1(4,4)&=&\frac{S^2_1(2,2)-S_1(1,2)}{2}-4-S_9(1,2)
\end{eqnarray*}
In the calculation of $S_1(2,4)$ and $S_1(4,4)$ 
on the basis of the relevant quadratic equation it is to consider that
$S_1(2,4)>S_1(4,4)$.

It remains only to notice that we need for the calculations in this step
the values $S_1(1,2)$, $S_1(2,2)$ as well as
$S_9(1,2)$ and $S_9(2,2)$. It is necessary to get them in step 5. \newline

{\bf Step 5.} We split at first $\;\;S_1=p_1+p_2+p_4+p_8+p_{16}+p_{32}+p_{64}+p_{128}$
into
$S_{1}(1,2)=p_1+p_4+p_{16}+p_{64}$ and $S_{1}(2,2)=p_2+p_8+p_{32}+p_{128}$.  

In this case we have in
(\ref{Zeile1})
only the product $p_1p_4$ before the middle and the product $p_1p_{16}$ in the  middle.
The pair $p_3$ belongs to $S_2$, $p_5$ belongs to $S_8$, and $p_{15}$ belongs to $S_9$,
so that $H_0=S_2+S_8+S_9$. Thus, for $S_1(1,2)$ and $S_1(2,2)$, the following holds:
\begin{eqnarray*}
S_1(1,2)+S_1(2,2)&=&S_1\\
S_1(1,2)\cdot S_1(2,2)&=&\frac{S^2_1-S_1}{2}-8-(S_2+S_8+S_9)
\end{eqnarray*}
The values $S_1(1,2)$ and $S_1(2,2)$ can be uniquely determined using the relevant 
quadratic equation with consideration that $S_1(1,2)>S_1(2,2)$.

We come from $S_1(1,2)$ and $S_1(2,2)$ to $S_9(1,2)$ and $S_9(2,2)$ 
via shift with the height $\textit{shift}_2=8$ for the block numbers:
\begin{eqnarray*}
S_9(1,2)+S_9(2,2)&=&S_9\\
S_9(1,2)\cdot S_9(2,2)&=&\frac{S^2_9-S_9}{2}-8-(S_1+S_{10}+S_{16})
\end{eqnarray*}
and it is to consider in the calculation that $S_9(2,2)>S_9(1,2)$.

The calculations in this step show that to determine $p_1$ we need 
the values $S_1$, $S_2$, $S_8$, $S_9$, $S_{10}$ and $S_{16}$.
This information can be used to avoid unneeded splits.
In step 4 it is only necessary to split $F(1,8)$, $F(2,8)$, and $F(8,8)$. 
To get $F(1,8)$, $F(2,8)$ and $F(8,8)$
we have to split only $F(1,4)$, $F(2,4)$, and $F(4,4)$ in step 3.
In step 2 it is necessary to calculate only these 3 values, 
we do not need the value $F(3,4)$. This completes the case $n=257$.

\section{Regular 65537-gon}

We need here 15 steps to obtain $p_1$. 
To reduce the number of unnecessary splits, we first consider the final 
steps involving splits for invariant blocks in reverse order.

In case $n=65537$, the calculation of the required coefficients $\mu_{k,2^m}$
and the comparison of the parts $F(k,2^m)$ and $S_k(j,2^m)$ were performed using 
a quite simple C program. The numerical accuracy is enough to guarantee 
that the assignments are correct. \newline

{\bf Step 15.} To get $p_1$, we have to split $S_1(1,8)=p_1+p_{256}$ 
into $p_1$ and $p_{256}$. In this case there are no products in
(\ref{Zeile1}),
and $H_3=0$. Thus it holds: 
\begin{eqnarray*}
p_1 + p_{256}&=& S_1(1,8)\\
p_1\cdot p_{256}&=&\frac{ S^2_1(1,8)-S_1(2,8)}{2}-2
\end{eqnarray*}
and $p_1$ and $p_{256}$ can be uniquely determined by solving 
the relevant quadratic equation and taking into account that $p_1>p_{256}$.
For this calculation, as we see,
we definitely need the values $S_1(1,8)$ and $S_1(2,8)$.\newline

{\bf Step 14.} In order to get $S_1(1,8)$, we have to split 
\[
S_1(1,4)=p_1+p_{16}+p_{256}+p_{4096}
\] 
into 
$S_1(1,8)=p_1+p_{256}$ and $S_1(5,8)=p_{16}+p_{4096}$.

In this case we have only the middle product $p_1p_{256}$ in
(\ref{Zeile1}).
The pair $p_{255}$ belongs to $S_{1025}(4,4)$, and $H_2=S_{1025}(4,4)$.

We come to the following equations for $S_1(1,8)$ and $S_1(5,8)\;$:
\begin{eqnarray*}
S_1(1,8)+S_1(5,8)&=&S_1(1,4)\\
S_1(1,8)\cdot S_1(5,8)&=&\frac{S^2_1(1,4)-S_1(2,4)}{2}-4-S_{1025}(4,4)
\end{eqnarray*}
and we can uniquely determine these values as solutions of the relevant quadratic equation
taking into account that $S_1(1,8)>S_1(5,8)$.

For $S_1(2,8)$ and $S_1(6,8)$  we have, compared to $S_1(1,8)$ and $S_1(5,8)$,
a shift for the pair numbers with the height $\textit{shift}_1=1$. Therefore we come
to the following equations for $S_1(2,8)$ and $S_1(6,8)\;$: 
\begin{eqnarray*}
S_1(2,8)+S_1(6,8)&=&S_1(2,4)\\
S_1(2,8)\cdot S_1(6,8)&=&\frac{S^2_1(2,4)-S_1(3,4)}{2}-4-S_{1025}(1,4)
\end{eqnarray*}
The values $S_1(2,8)$ and $S_1(6,8)$ can be uniquely determined using the corresponding 
quadratic equation and taking into account that $S_1(2,8)>S_1(6,8)$.

To perform the calculations in this step we obviously need the values
$S_1(1,4)$, $S_1(2,4)$, $S_1(3,4)$, $S_{1025}(1,4)$, and $S_{1025}(4,4)$. \newline

{\bf Step 13.} We obtain the values $S_1(1,4)$ and $S_1(3,4)$ if we 
split\newline 
$S_1(1,2)=p_1+p_4+p_{16}+p_{64}+p_{256}+p_{1024}+p_{4096}+p_{16384}$
\newline
into $S_1(1,4)=p_1+p_{16}+p_{256}+p_{4096}$ and $S_1(3,4)=p_4+p_{64}+p_{1024}+p_{16384}$.

Here we have in
(\ref{Zeile1})
the product $p_1p_{16}$ before the middle and $p_1p_{256}$ in the middle.
The pair $p_{15}$ belongs to $S_{1957}(2,2)$, $p_{17}$ belongs to $S_{1117}(2,2)$
and $p_{255}$ belongs to $S_{1025}(2,2)$, and hence
$H_1=S_{1957}(2,2)+S_{1117}(2,2)+S_{1025}(2,2)$.

Therefore we have:
\begin{eqnarray*}
S_1(1,4)+S_1(3,4)&=&S_1(1,2)\\
S_1(1,4)\cdot S_1(3,4)&=&\frac{S^2_1(1,2)-S_1(2,2)}{2}-8\\
&&\;-S_{1025}(2,2)-S_{1117}(2,2)-S_{1957}(2,2)
\end{eqnarray*}
The values $S_1(1,4)$ and $S_1(3,4)$ can be uniquely determined using the relevant
quadratic equation and taking into account that $S_1(1,4)>S_1(3,4)$.

We come from $S_1(1,4)$ and $S_1(3,4)$ to $S_1(2,4)$ and $S_1(4,4)$ via shift
with the height $\textit{shift}_1=1$ for the pair numbers.
For  $S_1(2,4)$ and $S_1(4,4)$ we obtain:
\begin{eqnarray*}
S_1(2,4)+S_1(4,4)&=&S_1(2,2)\\
S_1(2,4)\cdot S_1(4,4)&=&\frac{S^2_1(2,2)-S_1(1,2)}{2}-8\\
&&\;-S_{1025}(1,2)-S_{1117}(1,2)-S_{1957}(1,2)
\end{eqnarray*}
We can uniquely determine these values using the relevant 
quadratic equation and taking into account that $S_1(2,4)>S_1(4,4)$.

The value $S_{1025}(1,4)$ can be determined by splitting $S_{1025}(1,2)$ into\linebreak
$S_{1025}(1,4)$ and $S_{1025}(3,4)$. 
We come to $S_{1025}(1,4)$ and $S_{1025}(3,4)$ from $S_{1}(1,4)$ and $S_{1}(3,4)$ 
via shift with the height $\textit{shift}_2=1024$ for the block numbers.

Thus we obtain for the values $S_{1025}(1,4)$ and $S_{1025}(3,4)$:
\begin{eqnarray*}
S_{1025}(1,4)+S_{1025}(3,4)&=&S_{1025}(1,2)\\
S_{1025}(1,4)\cdot S_{1025}(3,4)&=&\frac{S^2_{1025}(1,2)-S_{1025}(2,2)}{2}-8\\
&&-S_1(1,2)-S_{93}(1,2)-S_{933}(1,2)
\end{eqnarray*}
and they can be uniquely determined using the relevant quadratic equation 
and taking into account that $S_{1025}(1,4)>S_{1025}(3,4)$. We need only $S_{1025}(1,4)$.

The value $S_{1025}(4,4)$ can be determined by splitting $S_{1025}(2,2)$ into\linebreak 
$S_{1025}(2,4)$ and $S_{1025}(4,4)$. 
Compared to $S_{1025}(1,4)$ and $S_{1025}(3,4)$ we have a shift
with the height $\textit{shift}_1=1$ for the pair numbers, and we obtain therefore:
\begin{eqnarray*}
S_{1025}(2,4)+S_{1025}(4,4)&=&S_{1025}(2,2)\\
S_{1025}(2,4)\cdot S_{1025}(4,4)&=&\frac{S^2_{1025}(2,2)-S_{1025}(1,2)}{2}-8\\
&&-S_1(2,2)-S_{93}(2,2)-S_{933}(2,2)
\end{eqnarray*}
These values can be uniquely determined using the relevant 
quadratic equation and taking into account that $S_{1025}(2,4)<S_{1025}(4,4)$. 
Here we need only the value $S_{1025}(4,4)$.

The calculations in this step show that to realize all of them, we need 
$S_1(1,2)$, $S_1(2,2)$, $S_{93}(1,2)$, $S_{93}(2,2)$,
$S_{255}(1,2)$, $S_{255}(2,2)$, $S_{933}(1,2)$, $S_{933}(2,2)$,
$S_{1025}(1,2)$, $S_{1025}(2,2)$, $S_{1117}(1,2)$, $S_{1117}(2,2)$,
$S_{1957}(1,2)$, and $S_{1957}(2,2)$.
To  obtain these values, we must split the blocks
$S_1$, $S_{93}$, $S_{225}$, $S_{933}$, $S_{1025}$, $S_{1117}$, and $S_{1957}$.
\newline

{\bf Step 12.} We split $S_1=p_1+p_2+p_4+\cdots+p_{16384}+p_{32768}$\newline
into $\; S_1(1,2)=p_1+p_4+p_{16}+p_{64}+p_{256}+p_{1024}+p_{4096}+p_{16384}$\newline
and
$ S_1(2,2)=p_2+p_8+p_{32}+p_{128}+p_{512}+p_{2048}-p_{8192}+p_{32768}\;$.

We have in 
(\ref{Zeile1})
the products $p_1p_{4},p_1p_{16}, p_1p_{64}$ before the middle
and $p_1p_{256}$ in the middle.
The pairs $p_3$, $p_5$, $p_{15}$, $p_{17}$, $p_{63}$, $p_{65}$, and $p_{255}$
belong to the invariant blocks 
$S_2$, $S_{1956}$, $S_{1957}$, $S_{1117}$, $S_{1266}$, $S_{1900}$, and $S_{1025}$.
Hence $H_0 = S_2 + S_{1025} + S_{1117} + S_{1266} + S_{1900} + S_{1956} + S_{1957}$.

For $S_1(1,2)$ and $S_1(2,2)$ we obtain:
\begin{eqnarray*}
S_1(1,2)+S_1(2,2)&=&S_1\\
S_1(1,2)\cdot S_1(2,2)&=&\frac{S^2_1-S_1}{2}-16-S_2-S_{1025}-S_{1117}-S_{1266}\\
&&\;\; -S_{1900}-S_{1956}-S_{1957}
\end{eqnarray*}
These values can be uniquely determined using the relevant quadratic equation 
and taking into account that $S_1(1,2)> S_1(2,2)$.

To get $S_{93}(1,2)$ and $S_{93}(2,2)$, we have to split $S_{93}$.
Compared to $S_{1}(1,2)$ and $S_{1}(2,2)$, we have a shift with the height 
$\textit{shift}_2=92$ for the block numbers.
We therefore have for $S_{93}(1,2)$ and $S_{93}(2,2)\;$:
\begin{eqnarray*}
S_{93}(1,2)+S_{93}(2,2)&=&S_{93}\\
S_{93}(1,2)\cdot S_{93}(2,2)&=&\frac{S^2_{93}-S_{93}}{2}-16-S_1-S_{94}-S_{1117}-S_{1209}\\
&&\;\; -S_{1358}-S_{1992}-S_{2048}
\end{eqnarray*}
We can uniquely determine these values by solving the relevant quadratic equation 
and taking into account that $S_{93}(1,2)< S_{93}(2,2)$.

For $S_{933}(1,2)$ and $S_{933}(2,2)$ we have, compared to $S_{1}(1,2)$ and $S_{1}(2,2)$,
the shift with the height $\textit{shift}_2=932$ for the block numbers. 
Therefore we obtain for $S_{933}(1,2)$ and $S_{933}(2,2)$:
\begin{eqnarray*}
S_{933}(1,2)+S_{933}(2,2)&=&S_{933}\\
S_{933}(1,2)\cdot S_{933}(2,2)&=&\frac{S^2_{933}-S_{933}}{2}-16-S_1-S_{150}-S_{784}\\
&&\;\; -S_{840}-S_{841}-S_{934}-S_{1957}
\end{eqnarray*}
These values can be uniquely determined using
the relevant quadratic equation and taking into account that $S_{933}(1,2)< S_{933}(2,2)$.

For $S_{1025}(1,2)$ and $S_{1025}(2,2)$ we have, compared to $S_{1}(1,2)$ and $S_{1}(2,2)$,
the shift with the height $\textit{shift}_2=1024$ for the block numbers. 
We get therefore for $S_{1025}(1,2)$ and $S_{1025}(2,2)\;$:
\begin{eqnarray*}
S_{1025}(1,2)+S_{1025}(2,2)&=&S_{1025}\\
S_{1025}(1,2)\cdot S_{1025}(2,2)&=&\frac{S^2_{1025}-S_{1025}}{2}-16-S_1-S_{93}-S_{242} \\
&&\;\;-S_{876}-S_{932}-S_{933}-S_{1026}
\end{eqnarray*}
We can uniquely determine these values using the relevant quadratic equation 
and taking into account that $S_{1025}(1,2)< S_{1025}(2,2)$.

For $S_{1117}(1,2)$ and $S_{1117}(2,2)$ we have, compared to $S_{1}(1,2)$ and $S_{1}(2,2)$,
the shift with the height $\textit{shift}_2=1116$ for the block numbers. 
We obtain for $S_{1117}(1,2)$ and $S_{1117}(2,2)\;$:
\begin{eqnarray*}
S_{1117}(1,2)+S_{1117}(2,2)&=&S_{1117}\\
S_{1117}(1,2)\cdot S_{1117}(2,2)&=&\frac{S^2_{1117}-S_{1117}}{2}-16 -S_{93}-S_{185}-S_{334} \\
&&\;\; -S_{968}-S_{1024}-S_{1025}-S_{1118}
\end{eqnarray*}
These values can be uniquely determined using the relevant quadratic equation 
and taking into account that $S_{1117}(1,2)< S_{1117}(2,2)$.

At last, for $S_{1957}(1,2)$ and $S_{1957}(2,2)$ we have, compared to 
$S_{1}(1,2)$ and $S_{1}(2,2)$,
the shift with the height $\textit{shift}_2=1956$ for the block numbers.
For $S_{1957}(1,2)$ and $S_{1957}(2,2)$ we obtain therefore:
\begin{eqnarray*}
S_{1957}(1,2)+S_{1957}(2,2)&=&S_{1957}\\
S_{1957}(1,2)\cdot S_{1957}(2,2)&=&\frac{S^2_{1957}-S_{1957}}{2}-16-S_{933}-S_{1025}-S_{1174}\\
&&\;\;-S_{1808}-S_{1864}-S_{1865}-S_{1958}
\end{eqnarray*}
These values can be uniquely determined using the relevant quadratic equation and
taking into account that $S_{1957}(1,2)< S_{1957}(2,2)$.

We have presented all splits in step 12 without gaps and can see that
we need for these splits the invariant blocks $S_j$ and $S_{j+1024}$  
for the following $18$ numbers $j$
\begin{eqnarray}\label{ZuTeilenF1024}
1,2,93,94,150,185,242,334,784,840,841,876,932,933,934,941,968,1024\;\;
\end{eqnarray}
$S_j$ and $S_{j+1ß24}$ are parts of $F(j,1024)$, and these 18 values $F(j,1024)$
we unconditionally have to split in step 11.

In step 11 we already split the values $F(j,1024)$. We consider this step in advance 
as we do not need any special calculation for these splits, and the results of these splits 
will also help us reduce the splittings in the other steps. \newline

{\bf Step 11.} We consider at first the split of $F(1,1024)=S_1+S_{1025}$
into $F(1,2048)=S_1$ and $F(1025,2048)=S_{1025}$.

The product $F(1,2048)\cdot F(1025,2048)=S_1\cdot S_{1025}$ is known,
(\ref{G1G1025}), 
and with the help of an obvious rearrangement we obtain for  $F(j,2048)$ and\linebreak 
$F(1024+j,2048)$, $1\le j\le 1024\;$: 
\begin{eqnarray}\label{F1F1025,2048}
&&F(j,2048)+ F(j+1024,2048)\;=\; F(j,1024)\nonumber\\
&&F(j,2048)\cdot F(j+1024,2048)\;=\;2F(j,1024)+F(\rho(23+j,1024),1024)\nonumber\\
&&\;\;+2F(\rho(154+j,1024),1024)+F(\rho(184+j,1024),1024)\nonumber\\
&&\;\;+F(\rho(308+j,1024),1024)+F(\rho(359+j,1024),1024)\nonumber\\
&&\;\;+F(\rho(530+j,1024),1024)+F(\rho(666+j,1024),1024)\\
&&\;\;+F(\rho(718+j,1024),1024)+F(\rho(733+j,1024),1024)\nonumber\\
&&\;\;+2F(\rho(777+j,1024),1024)+F(\rho(840+j,1024),1024)\nonumber\\
&&\;\;+F(\rho(945+j,1024),1024)\nonumber
\end{eqnarray}
We have seen in step 12 
that we have to split only 18 values $F(j,1024)$ with the numbers $j$ from the list
(\ref{ZuTeilenF1024}).
To uniquely determine the values\linebreak
$S_j=F(j,2048)$ and $S_{1024+j}=F(1024+j,2048)$, we need to know which of them is greater.
The following list contains only those of these $18$ numbers $j$ 
for which $F(j,2048)>F(1024+j,2048)\,$:
\[
1,\; 2,\;185,\;334,\;968,\;1024
\]
This information is doubtlessly sufficient to clearly determine  
$F(j,2048)$ and $F(1024+j,2048)$ with the help of the relevant quadratic equation. 

Now we want to analyze which of the values $F(j,1024)$ we really need to perform
the calculations in this step. To split $\;F(j,1024)\;$ for the 18 numbers $j$ from the list
(\ref{ZuTeilenF1024}),
we need the appropriate values $\; F(k,1024)$ to represent $F(j,2048)\cdot F(j+1024,2048)$.
The products $F(j,2048)\cdot F(j+1024,2048)$ are shown in 
(\ref{F1F1025,2048}),
and the corresponding numbers $k$ must be calculated for the 18 numbers $j$ from the list
(\ref{ZuTeilenF1024})
based on this representation.

So we see that we need, in total, 213 values $F(k,1024)$ to realize
the splittings in step 11. It is a lot of values, but clearly less than 1024 possible values.
The list of the 213 numbers $j$ (we use instead of $k$ our usual notation $j$)\linebreak
of the required values $F(j,1024)$ is the following: 
\begin{eqnarray*}
&&1, 2, 6, 14, 15, 23, 24, 25, 28, 36, 43, 58, 62, 63, 64, 68, 71, 87, 92, 93, 94, 98,\\
&&101, 106, 116, 117, 119, 124, 125, 128, 150, 154, 155, 156, 160, 163, 173, 175,\\
&&176, 184, 185, 186, 208, 211, 216, 217, 218, 225, 242, 247, 248, 252, 255,265,\\
&&267, 268, 269, 276, 277, 278, 290, 303, 304, 308, 309, 310, 334, 339, 346, 347,\\
&&357, 359, 360, 361, 369, 382, 396, 401, 402, 426, 438, 396, 401, 402, 426, 438,\\
&&439, 440, 447, 452, 453, 458, 474, 478, 482, 483, 488, 493,  509, 518, 530, 531,\\
&&532, 534, 535, 537, 544, 549, 550, 570, 574, 575, 576, 583, 585, 593, 594, 600,\\
&&601, 610, 623, 624, 626, 627, 628, 629, 635, 641, 642, 643, 650, 656, 657, 662,\\
&&666, 667, 668, 677, 680, 685, 686, 687, 692, 693, 694, 705, 715, 718, 719, 720,\\
&&721, 733, 734, 735, 748, 749, 750, 757, 759, 760, 761, 762, 772, 777, 778, 779,\\
&&784, 797, 807, 811, 812, 816, 826, 827, 840, 841, 842, 851, 853, 854, 855, 862,\\
&&863, 864, 868, 870, 871, 876, 883, 889, 899, 903, 908, 918, 927, 932, 933, 934,\\
&&938, 941, 945, 946, 947, 955, 956, 957, 960, 962, 964, 968, 975, 990, 991, 994,\\
&& 995, 1000, 1019, 1024
\end{eqnarray*}
For 32 of these 213 numbers, namely for the following numbers $j\;$:
\begin{eqnarray}\label{Nr_j}
&& 6, 23, 25, 58, 62, 63, 64, 71, 98, 116, 117, 150, 154,155, 156, 173, 175, 208,\;\;\nonumber\\
&& 247, 248, 265, 267, 304, 339, 359, 396, 426, 452, 478, 482, 483, 488\; , 
\end{eqnarray}
we need $F(j,1024)$ and $F(j+512, 1024)$.
For these numbers $j$ we obtain therefore these two necessary values if we split  
$F(j,512)$. Hence we need to do only $181$ splits in step 10. 
\newline

From here on we consider the steps in natural order.
Thereby we will use the already gained information about the required 
values $F(j,1024)$ to reduce the number of splits.
\newpage

{\bf Step 1.} We split $S$ into $F(1,2)$ and $F(2,2)$.

Here the sum and the product are known: $F(1,2)+F(2,2)=-1$,\linebreak
 $ F(1,2)\cdot F(2,2)=-(n-1)/4=-16384\;$. 
Therefore $F(1,2)$ and $F(2,2)$ can be uniquely determined using the quadratic equation
\begin{eqnarray*}\label{x^2+x-16384}
x^2+x-16384=0\; 
\end{eqnarray*}
and taking into account that $F(1,2)>F(2,2)$. \newline

{\bf Step 2.} For the split $F(j,2)$ into $F(j,4)$ and $F(2+j,4)$
the sum and the product are known:\newline
$F(j,4)+F(2+j,4)=F(j,2)$ and $F(j,4)\cdot F(2+j,4)=-(n-1)/16=-4096$.

Therefore $F(j,4)$ and $F(2+j,4)$, $1\le j\le 2$,
can be uniquely determined with the help of the quadratic equation
\[
x^2-F(j,2)x-4096=0\; ,
\]
considering that $F(j,4)<F(2+j,4)$, $1\le j\le 2$. \newline

{\bf Step 3.} First we look at the split of
\[
F(1,4)=S_1+S_5+S_9+S_{13}+\dots +S_{2041}+S_{2045}
\]
into $F(1,8)=S_1+S_9+\cdots +S_{2041}$ and $F(5,8)=S_5+S_{13}+\cdots +S_{2045}\;$.

The product $F(1,8)\cdot F(5,8)$
is a linear combination of the values $F(k,4)$ with integer coefficients
$\mu_{k,4}$, $1\le k\le 4$. These coefficients are 
\[
\mu_{1,4}=992,\; \mu_{2,4}=1040,\; \mu_{3,4}=1024,\; \mu_{4,4}=1040\; ,
\]
and it holds therefore
\begin{eqnarray*}
&&F(1,8)\cdot F(5,8)=992F(1,4)+1040 F(2,4)+1024 F(3,4)+1040 F(4,4)\\
&&\;\;=1040\cdot  S  -48F(1,4)-16F(3,4)\; =\; -1040 -48F(1,4)-16F(3,4)\; .
\end{eqnarray*}

For the parts $F(j,8)$ and $F(4+j,8)$, $1\le j\le 4$, we then obtain:
\begin{eqnarray*}
F(j,8)+ F(4+j,8)&=& F(j,4)\\
F(j,8)\cdot F(4+j,8)&=& -1040 -48F(j,4)-16F(\rho(2+j,4),4)
\end{eqnarray*}
and these values can be uniquely determined taking into account that
\[
F(1,8)<F(5,8),\; F(2,8)>F(6,8),\; F(3,8)>F(7,8),\; F(4,8)<F(8,8)\; .
\]

{\bf Step 4.} First we look at the split of
\[
F(1,8)=S_1+S_9+S_{17}+S_{25}+\cdots +S_{2033}+S_{2041}
\]
into $F(1,16)=S_1+S_{17}+\cdots +S_{2033}$ and $F(9,16)=S_9+S_{25}+\cdots +S_{2041}\;$.

The product $F(1,16)\cdot F(9,16)$
is a linear combination of the values $F(k,8)$ with integer coefficients
$\mu_{k,8}$, $1\le k\le 8$. The coefficients are
\begin{eqnarray*}
&\mu_{1,8}=284,\; \mu_{2,8}=237, \; \mu_{3,8}=272,\; \mu_{4,8}=237,\\
&\mu_{5,8}=256,\; \mu_{6,8}=269,\; \mu_{7,8}=256,\; \mu_{8,8}=237
\end{eqnarray*}
and we get therefore
\begin{eqnarray*}
&&F(1,16)\cdot F(9,16)=\sum_{k=1}^{8}\mu_{k,8}F(k,8)\\
&&\;\;\;=-237 +47 F(1,8)+35F(3,8)+19F(5,8)+32F(6,8)+19F(7,8)\; .
\end{eqnarray*}

For $F(j,16)$ and $F(8+j,16)$, $1\le j\le 8$, we obtain:
\begin{eqnarray*}
&&F(j,16)+ F(8+j,16)= F(j,8)\\
&&F(j,16)\cdot F(8+j,16)= -237+47 F(j,8)+35F(\rho(2+j,8),8)\\
&&\;\;\;+19F(\rho(4+j,8),8)+32F(\rho(5+j,8),8)+19F(\rho(6+j,8),8)
\end{eqnarray*}
and these values can be uniquely determined considering that
\begin{eqnarray*}
&&F(1,16)>F(9,16),\; F(2,16)>F(10,16),\; F(3,16)>F(11,16),\\
&& F(4,16)<F(12,16),\; F(5,16)>F(13,16),\; F(6,16)<F(14,16),\\
&&F(7,16)>F(15,16),\; F(8,16)<F(16,16).
\end{eqnarray*}

{\bf Step 5.} First we look at the split of
\[
F(1,16)=S_1+S_{17}+S_{33}+S_{49}+\cdots+S_{2017}+S_{2033}
\]
into $F(1,32)=S_1+S_{33}+\cdots +S_{2017}$ and $F(17,32)=S_{17}+S_{49}+\cdots +S_{2033}\,$.

The product $F(1,32)\cdot F(17,32)$
is a linear combination of the values $F(k,16)$ with integer coefficients
$\mu_{k,16}$, $1\le k\le 16$. These coefficients are
\begin{eqnarray*}
&&\mu_{1,16}=80,\; \mu_{2,16}=62, \; \mu_{3,16}=60,\; \mu_{4,16}=64,\; \mu_{5,16}=57,
\; \mu_{6,16}=60, \\
&&\mu_{7,16}=61,\; \mu_{8,16}=60,\;\mu_{9,16}=68,\; \mu_{10,16}=64,\;\mu_{11,16}=64,\; \mu_{12,16}=58,\\
&&\mu_{13,16}=65,\; \mu_{14,16}=70,\; \mu_{15,16}=61,\; \mu_{16,16}=70
\end{eqnarray*}
and therefore we obtain
\begin{eqnarray*}
&&F(1,32)\cdot F(17,32)=\sum_{k=1}^{16}\mu_{k,16}F(k,16)\\
&&=-60 +20 F(1,16)+2F(2,16)+4F(4,16)-3F(5,16)+F(7,16)\\
&&\;\;\;\;\; +8F(9,16)+4F(10,16)+4F(11,16)-2F(12,16)+5 F(13,16)\\
&&\;\;\;\;\; +10F(14,16)+F(15,16)+10F(16,16)\; .
\end{eqnarray*}
For $F(j,32)$ and $F(16+j,32)$, $1\le j\le 16$, we have therefore:
\begin{eqnarray*}
&&F(j,32)+ F(16+j,32)= F(j,16)\\
&&F(j,32)\cdot F(16+j,32)= -60 +20 F(j,16)+2F(\rho(1+j,16),16)\\
&&\;\; +4F(\rho(3+j,16),16)-3F(\rho(4+j,16),16) +F(\rho(6+j,16),16)\\
&&\;\; +8F(\rho(8+j,16),16)+4F(\rho(9+j,16),16)+4F(\rho(10+j,16),16)\\
&&\;\;-2F(\rho(11+j,16),16)+5 F(\rho(12+j,16),16)+10F(\rho(13+j,16),16)\\
&&\;\;+F(\rho(14+j,16),16)+10F(\rho(15+j,16),16)
\end{eqnarray*}
To determine $F(j,32)$ and $F(16+j,32)$, $1\le j\le 16$, 
using the relevant quadratic equation we have to know which of them is greater.
To summarize this information in short, we present in the following list 
only the numbers $j$, $1\le j\le 16$, for which we have $F(j,32)>F(16+j,32)\;$:
\[
1,\; 3,\; 4,\; 5,\; 8,\; 9,\; 10,\; 12,\; 13,\; 15  
\]
For numbers $\,j$, $\, 1\le j\le 16\,$, that are not present in this list 
it holds $F(j,32)<F (16+j,32)$. 
This information is obviously sufficient to uniquely determine the values
$F(j,32)$ and $F (16+j,32)$.
\newline

{\bf Step 6.} We consider at first the split of 
\[
F(1,32)=S_1+S_{33}+S_{65}+S_{97}+\cdots +S_{1985} +S_{2017}
\]
into $F(1,64)=S_1+S_{65}+\cdots +S_{1985}$ and $F(33,64)=S_{33}+S_{97}+\cdots +S_{2017}$

The product $F(1,64)\cdot F(33,64)$ is a linear combination of the values $F(k,32)$
with the  integer coefficients $\mu_{k,32}$, $1\le k\le 32$. These coefficients are
\begin{eqnarray*}
&&\mu_{1,32}=4,\, \mu_{2,32}=12,\, \mu_{3,32}=20,\, \mu_{4,32}=13,\, \mu_{5,32}=20,
\,\mu_{6,32}=18,\\
&&\mu_{7,32}=16,\, \mu_{8,32}=19,\, \mu_{9,32}=19,\, \mu_{10,32}=22,\,
\mu_{11,32}=12,\,\mu_{12,32}=22,\\
&&\mu_{13,32}=13,\, \mu_{14,32}=13,\, \mu_{15,32}=11,\,\mu_{16,32}=22,\, 
\mu_{17,32}=20,\,\mu_{18,32}=15,\\
&&\mu_{19,32}=25,\,\mu_{20,32}=12,\; \mu_{21,32}=16,\; \mu_{22,32}=12,\,\mu_{23,32}=16,\,
\mu_{24,32}=17,\\
&& \mu_{25,32}=29,\, \mu_{26,32}=16,\, \mu_{27,32}=7,\,\mu_{28,32}=17,\, \mu_{29,32}=13,
\ \mu_{30,32}=17,\\
&&\mu_{31,32}=13,\, \mu_{32,32}=11\; .
\end{eqnarray*}
The possible transformations of the representation of
$F(1,64)\cdot F(33,64)$ do not give significant simplifications.
Therefore we use the determined values $\mu_{k,32}$, $1\le k\le 32$, without transformation.

For $F(j,64)$ and $F(32+j,64)$, $1\le j\le 32$, we obtain:
\begin{eqnarray*}
&&F(j,64)+ F(32+j,64)= F(j,32)\\
&&\;\;F(j,64)\cdot F(32+j,64)=\sum_{k=1}^{32}\mu_{k,32}F(\rho(k+j-1,32),32)
\end{eqnarray*}
In the following list we present the numbers $j$, $1\le j\le 32$, for which 
we have $F(j,64)>F(32+j,64)$:
\[
1,\; 2,\; 4,\; 5,\; 9,\; 10,\; 11,\; 12,\; 21,\; 24,\; 28,\; 29,\; 31,\; 32\; .
\]
On the basis of this information $F(j,64)$ and $F(32+j,64)$, $1\le j\le 32$,
can be uniquely determined using the relevant quadratic equation. \newline

{\bf Step 7.} We consider at first the split of
\[
F(1,64)=S_1+S_{65}+S_{129}+S_{193}+\cdots +S_{1921} +S_{1985}
\]
into $\;F(1,128)=S_1+S_{129}+\cdots +S_{1921}$ and
$\;F(65,128)=S_{65}+S_{193}+\cdots +S_{1985}$.

The product $F(1,128)\cdot F(65,128)$
is a linear combination of the values $F(k,64)$ with integer coefficients
$\mu_{k,64}$, $1\le k\le 64$.
The list of the numbers $\;\mu_{k,64}$, $1\le k\le 64$, 
for this sum is very large, but these numbers are between $0$ and $10$.
In order to summarize this result we denote $\mathcal K(m,64)$, $1\le m\le 10$, 
the set of numbers $k$ for which applies $\mu_{k,64}=m$. 
A set  $\mathcal K(0,64)$ is not needed. It holds:
\begin{eqnarray*}
\mathcal K(1,64)&=&\{13,24,31,33,37,38\},\\
\mathcal K(2,64)&=&\{3,7,9,21,36,46,56,57\}\\
\mathcal K(3,64)&=&\{2,15,19,25,27,28,35,40,41,42,45,47,59,61,64\},\\
\mathcal K(4,64)&=&\{29,32,43,48,51,52,58,60\},\\
\mathcal K(5,64)&=&\{4,5,6,8,11,16,20,22,23,30,53,62\},\\
\mathcal K(6,64)&=&\{10,12,14,17,34,44,54\},\;\;\\
\mathcal K(7,64)&=&\{39,49,50\},\\
\mathcal K(8,64)&=&\{18,55,63\},\\
\mathcal K(9,64)&=& \emptyset ,\\
\mathcal K(10,64)&=&\{26\}\, .
\end{eqnarray*}

For $F(j,128)$ and $F(64+j,128)$, $1\le j\le 64$, we obtain:
\begin{eqnarray*}
F(j,128)+ F(64+j,128)&=& F(j,64)\\
F(j,128)\cdot F(64+j,128)&=&\sum_{m=1}^{10} m\left(\sum_{k\in \mathcal K(m,64)}
F(\rho(k+j-1,64),64)\right)\;\;
\end{eqnarray*}
In the following list we present the numbers $j$, $1\le j\le 64$, 
for which we have $F(j,128)>F(64+j,128)$:
\begin{eqnarray*}
&&1,\; 2,\; 3,\; 4,\; 6,\; 8,\; 10,\; 12,\; 16,\; 17,\; 18,\; 19,\; 20,\; 21,\;
24,\; 26,\; 27,\; 28,\; 29,\; 31,\\
&&32,\; 33,\; 34,\; 
36,\; 37,\; 39,\; 40,\; 42,\; 45,\; 46,\; 48,\; 50,\; 51,\; 57,\; 58,\; 59,\;
62,\; 63
\end{eqnarray*}
This information is sufficient to uniquely determine the values $F(j,128)$ and 
 $F(64+j,128)$ with the help of the relevant quadratic equation. \newline

{\bf Step 8.} We consider at first the split of 
\[
F(1,128)=S_1+S_{129}+S_{257}+\cdots +S_{1793}+S_{1921}
\]
into $\;F(1,256)=S_1+S_{257}+S_{513}+\cdots +S_{1793}$ and\newline 
$\;F(129,256)=S_{129}+S_{385}+S_{641}+\cdots +S_{1921}\;$.

The product $F(1,256)\cdot F(129,256)$ is a linear combination of the values
$F(k,128)$ with integer coefficients $\mu_{k,128}$, $1\le k\le 128$. 
These coefficients are between $0$ and $5$ and we denote $\mathcal K(m,128)$, $1\le m\le 5$, 
the set of numbers $k$ for which $\mu_{k,128}=m$. 
We have here:
\begin{eqnarray*}
\mathcal K(1,128)&=&\{2,4,5,7,8,9,15,16,17,21,23,26,27,31,36,38,46,48,49,52,\\
&& 57,59,61,62,81,83,87,90,91,96,99,100,111,112,116,117,\\
&& 119,120,124,125,126\},\\
\mathcal K(2,128)&=&\{1,34,35,37,39,41,42,43,60,64,68,71,74,75,77,80,84,88,\\
&& 89,95,98,102,104,105,109,110,118,122,128\},\\
\mathcal K(3,128)&=&\{11,30,101,106\},\\
\mathcal K(4,128)&=&\{66,107,115\},\\
\mathcal K(5,128)&=&\{58\}\; .
\end{eqnarray*}
For $F(j,256)$ and $F(128+j,256)$, $1\le j\le 128$, we obtain:
\begin{eqnarray*}
\!\!\!\!\!\!&&\!F(j,256)+ F(128+j,256)=F(j,128)\\
\!\!\!\!\!\!&&\!F(j,256)\cdot F(128+j,256)=\sum_{m=1}^{5}m\left (\sum_{k\in \mathcal K(m,128)}
F(\rho(k+j-1,128),128)\right)\;\;\;
\end{eqnarray*}
The list of numbers $j$ for which  $F(j,256)>F(128+j,256)$ is the following:
\begin{eqnarray*}
&&1,\, 2,\, 3,\, 4,\, 5,\, 6,\, 7,\, 9,\, 10,\, 11,\, 12,\, 14,\, 16,\, 20,\, 21,\, 24,\, 25,\, 26,\, 28,\, 29,\, 30,\,31,\\
&& 32,\, 33,\, 34,\, 35,\, 41,\, 44,\, 47,\, 49,\, 53,\, 55,\, 56,\, 60,\,
 61,\, 63,\, 65,\, 66,\, 67,\, 69,\,71,\\
&&72,\, 73,\, 74,\, 76,\, 77,\, 78,\, 82, \, 83,\, 85,\, 91,\, 92,\, 93,\, 98,\,
100,\, 102,\, 104,\, 105,\; 107\\
&&108,\; 109,\; 110,\;111,\; 113,\; 117,\; 120,\; 121,\; 124,\; 125,\, 126,\; 128\;. 
\end{eqnarray*}
This information is sufficient to uniquely determine $F(j,256)$ and\linebreak
$F(128+j,256)$, $1\le j\le 128$, using the relevant quadratic equation. \newline

{\bf Step 9.} At first we consider the split of
\[
F(1,256)=S_1+S_{257}+S_{513}+\cdots +S_{1537}+S_{1793}
\]
into $F(1,512)=S_1+S_{513}+S_{1025}+S_{1537}$ and\newline
$F(257,512)=S_{257}+S_{769}+S_{1281}+S_{1793}\,$.

The product $F(1,512)\cdot F(257,512)$ 
is a linear combination of the values $F(k,256)$ with integer coefficients
$\mu_{k,256}$, $1\le k\le 256$. The coefficients $\mu_{k,256}$ are between $0$ and $3$.
Analogously to the previous steps we denote for $1\le m\le 3$ the set 
$\mathcal K(m, 256)$ of numbers $k$ for which $\mu_{k,256}=m$. Here we have:
\begin{eqnarray*}
\mathcal K(1,256)&=&\{5,8,15,20,30,31,34,38,40,42,44,45,51,52,54,57,60,62,\\
&&66,69,71,79,80,82,85,89,90,107,110,113,118,125,129,136,\\
&&143,147,174,176,187,188,189,196,201,213,220,232,234,\\
&&244,251,253,254\},\\
\mathcal K(2,256)&=&\{4,29,157,186,246\},\\
\mathcal K(3,256)&=&\{83\}\; .
\end{eqnarray*}
For $F(j,512)$ and $F(256+j,512)$, $1\le j\le 256$, we obtain:
\begin{eqnarray*}
\!\!\!\!\!\!\;\;&&F(j,512)+ F(256+j,512)= F(j,256)\\
\!\!\!\!\!\!\;\;&&F(j,512)\cdot F(256+j,512)=\sum_{m=1}^{3}m\left(\sum_{k\in \mathcal K(m,256)}
F(\rho(k+j-1,256),256)\right),\;\;\;\;\;
\end{eqnarray*}
and the list of the numbers $j$ for which we have $F(j,512)>F(256+j,512)$ is the following:
\begin{eqnarray*}
&&1, 2, 3, 5, 6, 13, 14, 15, 16, 18, 19, 20, 21, 22, 23, 24, 25, 27, 30, 31, 32,
37, 38,\\
&&39, 42, 44, 45, 46, 47, 48, 49, 52, 54, 55, 58, 59, 62, 63, 65, 66, 76,
77, 78, 79,\\
&&82, 84, 87, 89, 90, 91, 93, 94, 96, 97, 98, 102, 103, 107, 110, 115, 117, 118, 119,\\
&&123,127, 133, 134, 136, 139, 144, 145, 147, 148, 149, 151, 152, 155, 156, 157,\\
&&160, 161, 162, 166, 168, 170, 171, 180, 182, 185, 186, 187, 188, 189, 191, 194,\\
&&195, 197, 200, 203, 206, 207, 209, 210, 213, 215, 217, 220, 221, 222, 223, 225,\\
&&230, 232, 234, 235, 236, 237, 238, 240, 241, 242, 244, 245, 254, 255, 256\, .
\end{eqnarray*}
This information is sufficient to uniquely determine the values $F(j,512)$ and 
 $F(256+j,512)$, $1\le j\le 256$, using the relevant quadratic equation. \newline

{\bf Step 10.} We consider at first the split of 
\[
F(1,512)=S_1+S_{513}+S_{1025}+S_{1537}
\]
into $F(1,1024)=S_1+S_{1025}$ and $F(513,1024)=S_{513}+S_{1537}$.

The Product $F(1,1024)\cdot F(513,1024)$ is a linear combination of the values 
$F(k,512)$ with integer coefficients $\mu_{k,512}$, $1\le k\le 512$,
These coefficients are between $0$ and $2$,
and we denote similarly to the previous for $1\le m\le 2$ the set $\mathcal K(m,512)$
of the numbers $k$ for which $\mu_{k,512}=m$. Here we have:
\begin{eqnarray*}
\mathcal K(1,512)&=&\{41,49,81,92,106,109,114,211,226,233,269,275,278,\\
&&281,303,349,379,390,431,465\},\\
\mathcal K(2,512)&=&\{68,86,88,135,175,451\}\; .
\end{eqnarray*}
For $F(j,1024)$ and $F(512+j,1024)$, $1\le j\le 512$, we obtain therefore:
\begin{eqnarray*}
\!\!\!\!\!\!\!\!&& F(j,1024)+ F(152+j,1024)=F(j,512)\\
\!\!\!\!\!\!\!\!&& F(j,1024)\cdot F(512+j,1024)=\sum_{m=1}^{2}m\left( \sum_{k\in \mathcal K(m,512)}
F(\rho(k+j-1,512),512)\right)
\end{eqnarray*}
In step 11 we have seen that only 213 values $F(j,1024)$ are unconditionally required 
and therefore we have to split 181 values $F(j,512)$. 
The list of the numbers $j$ for which we have to split $\; F(j,512)\;$
into $ \;F(j,1024)\;$ and $F(512+j,1024)$ is the following: 
\begin{eqnarray*}
&&1, 2, 6, 14, 15, 18, 19, 20, 22, 23, 24, 25, 28, 32, 36, 37, 38, 43, 58, 62, 63, 64,\\
&&68, 71, 73, 81, 82, 87, 88, 89, 92, 93, 94, 98, 101, 106, 111, 112, 114, 115, 116,\\
&&117, 119, 123, 124, 125, 128, 129, 130, 131, 138, 144, 145, 150, 154, 155, 156,\\
&&160, 163, 165, 168, 173, 174, 175, 176, 180, 181, 182, 184, 185, 186, 193, 203,\\
&&206, 207, 208, 209, 211, 216, 217, 218, 221, 222, 223, 225, 236, 237, 238, 242,\\
&&245, 247, 248, 249, 250, 252, 255, 260, 265, 266, 267, 268, 269, 272, 276, 277,\\
&&278, 285, 290, 295, 299, 300, 303, 304, 308, 309, 310, 314, 315, 328, 329, 330,\\
&&334, 339, 341, 342, 343, 346, 347, 350, 351, 352, 356, 357, 358, 359, 360, 361,\\
&&364, 369, 371, 377, 382, 387, 391, 396, 401, 402, 406, 415, 420, 421, 422, 426,\\
&&429, 433, 434, 435, 438, 439, 440, 443, 444, 445, 447, 448, 450, 452, 453, 456, \\
&&458, 463, 474, 478, 479, 482, 483, 488, 493, 507, 509, 512\; .
\end{eqnarray*}
In the following list we present the part of these 181 numbers $j$ 
for which $F(j,1024)>F(512+j,1024)$ holds:
\begin{eqnarray*}
&&1, 2, 6, 18, 22, 24, 25, 37, 62,  68, 73, 93, 94, 98, 101, 106, 112, 116, 124, 125,\\
&&131, 138, 144, 150, 154, 155, 156, 163, 165, 168, 175, 176, 180, 182, 184, 
185,\\
&&186, 206, 209, 211, 216, 218, 221, 225, 236, 238, 242, 260, 265, 268, 269, 
272,\\
&&276, 277, 295, 299, 300, 304, 309, 310, 314, 315, 328, 330, 334, 342, 350, 
352,\\
&&356, 357, 359, 360, 369, 371, 377, 382, 387, 396, 406, 415, 426, 429, 433, 
438,\\
&&440, 444, 447, 448, 452, 463, 474, 479, 482, 483, 512\; .
\end{eqnarray*}
This information is sufficient to uniquely determine the necessarily required
values 
$F(j,1024)$ and $F(512+j,1024)$. In these 181 splits of $F(j,512)$ we need both parts 
$F(j,1024)$ and $F(512+j,1024)$ only for 32 numbers $j$.
These 32 numbers $j$ we have seen in step 11, 
(\ref{Nr_j}). 

The calculations in step 10 do not indicate that the splittings in step 9
can be substantially reduced. Hence we have to do all splits in
the steps 1-9. This concludes the case $n=65537$. 


\newpage
\section{References}
\;\;\;\,
1. C. F. Gauss, Werke, Bd. 1, Disquisitiones arithmeticae, § 365, Königliche 
Gesellschaft der Wissenschaften zu Göttingen, 1863.

2. C. F. Gauss, Disquisitiones Arithmeticae, English translation by Arthur A. Clarke, 
New Haven, CT: Yale University Press, 1966.

3. M. G. Paucker, Geometrische Verzeichnung des regelmäßigen Siebzehn-Ecks 
und Zweyhundertsiebenundfunfzig-Ecks in den Kreis, Jahresverhandlung-en 
der Kurländischen Gesellschaft für Literatur und Kunst, 2, 1822, 160–219.

4. C. F. Gauss, Göttingische Gelehrte Anzeigen. Band 87, Nr. 203,\linebreak
19. Dezember 1825, 2025–2027.

5. F. J. Richelot, De resolutione algebraica aequationis $X^{257}=1$, sive de divisione 
circuli per bisectionem anguli septies repetitam in partes 257 inter se aequales commentatio 
coronata, J. reine angew. Math., 9, 1832, 1-26, 146-161, 209-230, and 337-358.

6. P. L. Wantzel, Recherches sur le moyens de reconnaıtre si un Problème de
Géométrie peut se résoudre avec la règle et le compas. Journal de Mathéma-tiques pures 
et appliquées 2, 1837, 366–372.

7. H. W. Richmond, A Construction for a regular polygon of seventeen sides, 
Quarterly Journal of Pure and Applied Mathematics, 26, 1893,\linebreak 206-207.

8. J. Hermes, Über die Teilung des Kreises in 65537 gleiche Teile, Nachrich-ten von der 
Gesellschaft der Wissenschaften zu Göttingen, Mathematisch-Physikalische Klasse, 
1894, 170–186, http://eudml.org/doc/58274

9. F. Klein, Vorträge über ausgewählte Fragen der Elementargeometrie, Leipzig, 1895.

10. H. S. Carslaw, Gauss's Theorem on the Regular Polygons which can be constructed 
by Euclid's Method, Cambridge University Press, 1910.

11. D. W. DeTemple, Carlyle Circles and the Lemoine Simplicity of Polygonal Constructions,
 Amer. Math. Monthly, 98, 1991, 97-108.
 
12. M. Trott, "$\cos(2\pi/257)$ à la Gauss." Mathematica Educ. Res., 4, 1995, 31-36.

13. C. Gottlieb, The Simple and Straightforward Construction of the Regular 257-gon, 
Mathematical Intelligencer, Vol. 21, No. 1, 1999, 31–37.

14. C. Goldstein, N. Schappacher, J. Schwermer,  
The Shaping of Arithmetic after C. F. Gauss's Disquisitiones Arithmeticae, 
Springer-Verlag, 2007.

15. H. M. Edwards, The construction of solvable polynomials, Bull. Amer. Math. Soc., 46,
2009, 397-411.

16. A. Skopenkov, Some More Proofs from the Book: Solvability and Insolvability of Equations 
in Radicals, 2014, https://arxiv.org/abs/0804.4357
\end{document}